\documentclass{article}[12pt]

\usepackage{epsfig}
\usepackage{amssymb}
\usepackage{amsmath,amsthm,amsfonts}
\usepackage{mathrsfs}
\usepackage{subfigure}
\usepackage{color}

\def\argmax{\mathop{\rm arg\,max}}

\def\argmin{\mathop{\rm arg\,min}}
\def\Real{\hbox{I\kern-.1667em\hbox{R}}}

\usepackage[sectionbib]{natbib}
\newcommand{\ldX}{X_{\bf h}}
\newcommand{\rX}{{\bf X}}
\newcommand{\h}{({\bf h})}

\newcommand{\blindcol}{black}
\renewcommand{\Re}{\Real}

\newcommand{\onecolwidth}{3.65in}
\newcommand{\twowidth}{2.8in}

\newcommand{\bb}{\mbox{\boldmath $\beta$}}
\newcommand{\btheta}{\mbox{\boldmath $\theta$}}

\newcommand{\XsubD}{X_{\mbox{\tiny D}}}

\newcommand{\supm}{\mbox{\tiny $\left(-1\right)$}}
\newcommand{\hmone}{{\bf h}^{\mbox{\tiny $\left(-1\right)$}}}

\theoremstyle{remark}
\newtheorem{thm}{Theorem} %insert [section] for 2 digit num
\newtheorem{prop}{Proposition}
\newtheorem{lem}{Lemma}
\newtheorem{coro}{Corollary}
\newtheorem{rem}{Remark}

\newcommand{\ass}[1]%
        {\renewcommand{\theenumi}{(#1)}%
        \item}

\title{Selecting Local Models in Multiple Regression\\ by Maximizing Power}
\author{\textcolor{\blindcol}{Chad M. Schafer and Kjell A. Doksum}\thanks{
\textcolor{\blindcol}{Chad M. Schafer is Visiting
Assistant Professor, Department of Statistics, Carnegie Mellon University,
Pittsburgh, PA 15213 (email: cschafer@stat.cmu.edu). 
Kjell A. Doksum is
Professor, Department of Statistics, University of Wisconsin, Madison,
WI 53607 (email: doksum@stat.wisc.edu). 
This research was
partially supported by NSF grants DMS-0505651, DMS-0604931, and DMS-0240019.
The authors are grateful to Alex Samarov for many helpful comments.}}}
\begin{document}

%\linenumbers
%----------------------------------------------------------------------------------------
% THE ABSTRACT
%----------------------------------------------------------------------------------------

\maketitle
\begin{abstract}
This paper considers multiple regression procedures for analyzing the relationship 
between a response variable and a vector of $d$ covariates in a nonparametric
setting where both tuning parameters and the number of covariates need to be selected.
We introduce an approach 
which handles the dilemma that with high dimensional data the sparsity of data in 
regions of the sample space makes estimation of nonparametric curves and surfaces virtually impossible.
This is accomplished by abandoning the goal of trying to estimate true underlying curves 
and instead estimating measures of dependence that can determine important relationships 
between variables. These dependence measures are based on local parametric fits on subsets 
of the covariate space that vary in both dimension and size within each dimension.
The subset which maximizes a signal to noise ratio is chosen, where the signal is a local 
estimate of a dependence parameter which depends on the subset dimension and size, and 
the noise is an estimate of the standard error (SE) of the estimated signal. This 
approach of choosing the window size to maximize a signal to noise ratio lifts the 
curse of dimensionality because for regions with sparsity of data the SE is very large.  
It corresponds to asymptotically maximizing the probability of correctly finding 
non-spurious relationships between covariates and a response or, more precisely, maximizing 
asymptotic power among a class of asymptotic level $\alpha$ $t$-tests indexed by subsets
of the covariate space. Subsets that achieve this goal are called features. 
We investigate the properties of specific procedures based on the preceding ideas 
using asymptotic theory and Monte Carlo simulations and find that within a 
selected dimension, the volume of the optimally selected subset does not tend to zero 
as $n \rightarrow \infty$ unless the volume of the subset of the covariate space 
where the response depends on the covariate vector tends to zero.  

\vspace{.1in}
\end{abstract}

%----------------------------------------------------------------------------------------
% THE INTRODUCTION
%----------------------------------------------------------------------------------------

\large
\section{Introduction}

In this paper, we analyze the relationship between a response variable $Y$ 
and a vector of covariates ${\bf X} = (X_1, X_2, \ldots , X_d)^T$ by using ``signal to noise'' 
descriptive measures of dependencies between $Y$ and ${\bf X}$.  This approach addresses a
dilemma of nonparametric statistical analysis (the curse of dimensionality), that with 
high dimensional data, sparsity of data in large regions of the sample space makes estimation 
of nonparametric curves and surfaces virtually impossible.  On the other hand, if the goal is
to find regions of strong dependence between variables, parametric methods may provide evidence for 
such relationships. This paper develops hybrid methods between parametric and nonparametric procedures 
that are designed to lift the curse of dimensionality by uniting the goal of analyzing dependencies 
between $Y$ and the $X$'s with the goal of finding a good model.

These methods are based on measures of dependencies between variables rather than on 
the estimation of some ``true'' curve or surface.  Such measures can, for instance, be based 
on the currently available procedures that are defined in terms of tuning parameters.  However, 
instead of asking what value of the tuning parameter will bring us closest to the ``truth'', 
we ask what value of the tuning parameter will give us the best chance of finding a relationship 
between the variables, if it exists.  For instance, consider the case where the tuning parameter is the window size 
${\bf h} = (h_1, \ldots , h_d)^T$ of local regions in $\Real^d$ within which we do local parametric fits 
of $Y$ to ${\bf X}$.  
Our approach consists of finding the ``best'' window size ${\bf h}$  by maximizing 
a local signal to noise ratio where the signal is an estimate of a measure of dependence and 
the noise is the estimated standard error (SE) of that estimate.
By dividing by the noise we lift the curse of dimensionality because the noise will be very 
large if the tuning parameter is such that the local region in the sample space is nearly empty.

Here we focus on the problem of exploring the relationship between the response $Y$ and
a {\em covariate of interest} $X_1$ while controlling for covariates $X_2, X_3, \ldots,X_d$.
%{\bf FIX THIS PASSAGE: We will consider both the case of partial regression of $X_1$ on the
%residuals of the fit of $Y$ on $X_2, X_3, \ldots, X_d$ and that of full regression of
%$Y$ on all $d$ covariates simultaneously. Which of these two approaches is a better choice
%will depend on the inference goals.}
Applied research abounds with situations where one is interested in
the relationship between two variables
while controlling for many other factors. The usual approach to dealing with the curse of
dimensionality is to assume either 
a simple parametric relationship (usually linear) or a nonparametric
but additive relationship between the covariates
and response. Our approach allows for nonparametric modeling of the interaction effects
and hence can find more subtle patterns in the data.

Our analysis methods begin with a {\em reduction step}, a procedure to analyze the multivariate
data in a family of subsets of the covariate space, subsets that vary in both
dimension and size within each dimension.
This step is motivated by the question ``Over which subsets of the covariates do these 
data provide evidence of a significant relationship between the response and the covariate of
interest?''
Once these subsets (called {\em features}) are extracted from the data, the analysis can go
in different directions depending on the objective. Inspection of the features is itself
a useful exploratory data analysis tool.
Scale space analysis is possible by varying the minimum feature size.
%The procedure can be used as a basis for a test of the 
%hypothesis that $X_1$ is an important predictor variable.

%WHERE THE INTRO ORIGINALLY ENDED
%\section{Local Testing and Asymptotic Power}\label{loctest}

The idea of using the signal to noise to choose a good procedure is
motivated by the result (\citet*{Pitman1948}, \citet*{Pitman1979},
\citet*{Serfling1980}, \citet*{Lehmann1999}) that the
asymptotic power of asymptotically normal test statistics $T_n$
for testing $H_0\!: \theta = \theta_0$ versus $H_1\!: \theta > \theta_0$
when $\theta$ is contiguous to $\theta_0$ can be compared by considering
their {\em efficacies}, defined as
\begin{equation}
   {\rm EFF}_0(T_n) = \frac{
   \frac{\partial}{\partial\theta} E_\theta(T_n) |_{\theta = \theta_0}}
   {{\rm SD}_{\theta_0}(T_n)},
\end{equation}
or the pre-limit version
\begin{equation}
   {\rm EFF}_\theta(T_n) = \frac{E_\theta(T_n) - E_{\theta_0}(T_n)}
   {{\rm SD}_{\theta}(T_n)}.
\end{equation}
See \citet*{Doksum2006}, who considered the single covariate case and selected
subsets using estimates of the nonparametric efficacy
\begin{equation}
   {\rm EFF}(T_n) = \frac{E_P(T_n)-E_{H_0}(T_n)}{{\rm SD}(T_n)}
\end{equation}
where $H_0$ is a specified independence property, $P$ is the probability distribution
of $({\bf X}, Y)$, and ${\rm SD}(T_n)$ is either ${\rm SD}_P(T_n)$ or
${\rm SD}_{H_0}(T_n)$. Let ${\rm SE}(T_n)$ denote a consistent estimator for
${\rm SD}(T_n)$ in the sense that ${\rm SE}(T_n)/{\rm SD}(T_n) \stackrel{P}{\longrightarrow}
1$.
We refer to estimates 
\begin{equation}
    \widehat{\rm EFF} = \frac{T_n - \widehat E_{H_0}(T_n)}{{\rm SE}(T_n)}
\end{equation}
of ${\rm EFF}(T_n)$ as a {\it signal to noise ratio} or {\it $t$-statistic}.
If $T_n \stackrel{P}{\longrightarrow} T(P)$ for some functional
$T(P) \equiv \theta$, then $\widehat{\rm EFF}$ resembles the familiar Wald test statistic
\begin{equation}
   \widehat t_W = \frac{\widehat \theta_n - \theta_0}{{\rm SE}(\widehat \theta_n)}.
\end{equation}

By Slutsky's Theorem, $\widehat {\rm EFF} \stackrel{\cal L}{\longrightarrow}
N(0,1)$ in general, as is the case of $\widehat t_W$ for parametric models.
We could refer to ``$\widehat {\rm EFF} \stackrel{\cal L}{\longrightarrow} N(0,1)$''
as the {\it Wald-Slutsky phenomenon}.
In some frameworks, if we consider $(\widehat {\rm EFF})^2$, this is related
to the Wilks phenomenon 
(\citet*{Fan2001}). 
%(Fan, Zhang, and Zhang (2001)).
By not squaring $\widehat {\rm EFF}$ we have a much simpler bias problem
%at the price of being limited to real-valued $T_n$. This simplification
that
makes it possible to address new optimality questions. 

Most work in the nonparametric testing literature starts with a class of tests
depending on a bandwidth ${\bf h}$ (or other tuning parameter) which tends to zero and
then finds the fastest possible rate that an alternative can converge to the null
hypothesis and still be detectable by one of the tests in the given class.
These rates and alternatives depend on ${\bf h}$ and the sample size $n$. See, for
instance, the development in Section 3.2 of \citet*{Fan2001}.
We consider alternatives {\it not} depending on ${\bf h}$ and ask what ${\bf h}$ maximizes the
probability of detecting the alternative. We find that
this ${\bf h}$ does not tend to ${\bf 0}$ unless the volume of the set on which $E(Y|{\bf X})$ is
nonconstant tends to zero as $n \rightarrow \infty$.

This paper is organized as follows.
%Section \ref{loctest} motivates the maximization of $t$-statistics in order
%to find significant features in the data.
Section \ref{loceff} defines and motivates our signal to noise criterion
and Section \ref{shrink} develops asymptotic properties.
Section \ref{varsel} describes how we incorporate variable selection into
the procedure for bandwidth selection.
Section \ref{critvals} explains how critical values can be approximated via
simulating the distribution of the criterion under random permutations of 
the data.
Section \ref{mc} shows results from Monte Carlo simulations and analysis of
real data.

\section{Local Efficacy in Nonparametric Regression}\label{loceff}

We consider $({\bf X}_1, Y_1), ({\bf X}_2, Y_2), \ldots, ({\bf X}_n, Y_n)$ i.i.d.\ as 
$({\bf X}, Y) \sim P$, ${\bf X} \in \Real^d$, $Y \in \Real$, and write
$Y = \mu({\bf X}) + \epsilon$ where $\mu({\bf X}) \equiv E(Y | {\bf X})$ is
assumed to exist and where $\epsilon = Y - \mu({\bf X})$. 
%We will at first focus on finding dependencies between $Y$ and ${\bf X}$ for 
%${\bf X}$ in a neighborhood of a given covariate vector ${\bf x}_0$
%from a unit (patient, component, DNA sequence) of interest. We want to know if a perturbation
%of ${\bf x}$ will affect the mean response for units with covariate values near
%${\bf x}_0$.
%Next we focus on the covariates one at a time, say $X_1$, and ask:
We will focus on finding dependencies between $Y$ and the covariate
of interest $X_1$ for 
${\bf X}$ in a neighborhood of a given covariate vector ${\bf x}_0$
from a targeted unit (e.g. patient, component, DNA sequence). 
We want to know if a perturbation
of $x_1$ will affect the mean response for units with covariate values near
${\bf x}_0$.
%Next we focus on the covariates one at a time, say $X_1$, and ask:
%We focus on one covariate (say $X_1$), 
%called the {\em covariate of interest}, and ask:
%Is there a relationship between $Y$ and $X_1$ at $x_1 = x_{01}$? 
Formally, we test
$H_0^{(1)}\!\!:$``$X_1$ is independent of $X_2, X_3, \ldots, X_d, Y$'' versus
$H_1^{(1)}\!\!:$``$\mu({\bf x})$ is not constant as a function of $x_1 \in \Real$.''

Our test statistics will be based on linear fits for ${\bf x}$ restricted to subregions of
the sample space. The subregion that best highlights the dependencies between $Y$ and covariates
will be determined by maximizing signal to noise ratios. We illustrate such procedures and their
properties by first considering
a linear model fit locally over the neighborhood $N_{\bf h}({\bf x}_0)$ 
of ${\bf x}_0$ where ${\bf h} = (h_1, h_2, \ldots, h_d)$ are bandwidths in each of the 
$d$ dimensions. The true 
relationship between $Y$ and ${\bf X}$ is such that
   $E[Y | \:{\bf X} = {\bf x}] \equiv \mu({\bf x})$ 
and ${\rm Var}(Y | \:{\bf X} = {\bf x}) \equiv \sigma^2$ are unknown.
The local linear model gives fitted values as a function of ${\bf x} \in N_{\bf h}({\bf x}_0)$, 
denoted 
\begin{equation}
   \widehat \mu_L({\bf x}) \equiv \widehat \beta_0 + \sum_{j=1}^d \widehat \beta_j \!\left(x_j - x_{0j}\right),
\end{equation}
where $\widehat \beta_j = \widehat \beta_j({\bf h})$ are sample coefficients depending on ${\bf x}_0$ and ${\bf h}$.
Thus, for coefficients $\beta_j= \beta_j({\bf h})$ depending on ${\bf x}_0$ and ${\bf h}$,
\begin{equation}
   \mu_{{\bf h}}\!\left({\bf x}\right) \equiv 
   E[\widehat \mu_L\!\left({\bf x}\right)] = \beta_0 + 
   \sum_{j=1}^d \beta_j \!\left(x_j-x_{0j}\right).
\end{equation}
We let $\widehat \beta_1({\bf h})$ be the test statistic $T_n$ for testing $H_0^{(1)}$ versus
$H_1^{(1)}$ and develop some of the properties of the efficacy and estimated efficacy for this
test statistic. Let ${\cal V}({\bf h}) \equiv P({\bf X} \in N_{\bf h}({\bf x}_0))$ and assume
throughout that $0 < {\cal V}({\bf h}) \leq 1$.

\subsection{The Signal\label{signal}}

Formally, the {\it signal}\/ is $E(\widehat \beta_1({\bf h})) = \beta_1({\bf h})$
where $\bb({\bf h}) = (\beta_0({\bf h}), \beta_1({\bf h}), \ldots, \beta_d({\bf h}))^T$ is the 
coefficient vector of the best local linear fit to $Y$. That is
\begin{equation}
   \bb({\bf h}) \equiv \argmin \{ E^{\bf h}[Y - (a + {\bf b}^T{\bf X})]^2\!: a \in \Real, \;{\bf b} \in \Real^d\}
   \label{optprob}
\end{equation}
where $E^{\bf h}$ is expected value for the conditional distribution $P^{\bf h}$
of $({\bf X}, Y)$ given
${\bf X} \in N_{\bf h}({\bf x}_0)$.
With ${\rm Var}^{\bf h}$ and ${\rm Cov}^{\bf h}$ also being conditional,
\begin{equation}
   \left( \beta_1({\bf h}), \ldots, \beta_d({\bf h})\right)^T = \Sigma_{\bf X}^{-1}({\bf h})
   \Sigma_{{\bf X}Y}({\bf h})
\end{equation}
where $\Sigma_{\bf X}({\bf h}) \equiv {\rm Var}^{\bf h}({\bf X})$ and
$\Sigma_{{\bf X}Y}({\bf h}) \equiv {\rm Cov}^{\bf h}({\bf X},Y)$.

Similarly, $\widehat \bb({\bf h})$ is the local linear least squares estimator
\begin{equation}
   \widehat \bb({\bf h}) = \left(\ldX^T \ldX\right)^{-1} \!\ldX^T {\bf Y}_{\bf h}
\end{equation}
where $\ldX$ is the design matrix $(X_{ij}), X_{i0} = 1, \:0 \leq j \leq d,\:
i \in {\cal I}$ with ${\cal I} = \{k\!: {\bf X}_k \in N_{\bf h}({\bf x}_0)\}$, and
${\bf Y}_{\bf h} = \{Y_i \!: i \in {\cal I}\}$.
It follows that conditionally given ${\mathbb X} \equiv ({\bf X}_1, \ldots, {\bf X}_n)$,
\begin{equation}
   {\rm Var}(\widehat \bb({\bf h})\:|\:{\mathbb X}) = \sigma^2 (\ldX^T \ldX)^{-1}.
\end{equation}
By writing $\ldX^T \ldX = \XsubD^T W_{\bf h} \XsubD$
where $W_{\bf h} \equiv {\rm diag}({\bf 1}[{\bf X}_i \in N_{\bf h}({\bf x}_0)]),
\: 0 \leq i \leq n$, with ${\bf 1}[{\bf X}_0 \in N_{\bf h}({\bf x}_0)] \equiv 1$,
and $\XsubD$ is the full design matrix $(X_{ij})_{n \times (d+1)}$,
we see that by the law of large numbers,
%with ${\cal V}({\bf h}) \equiv P({\bf X} \in N_{\bf h}({\bf x}))$,
\begin{eqnarray}
   n^{-1} \ldX^T \ldX \stackrel{a.s.}{\longrightarrow} 
   E(\XsubD \XsubD^T {\bf 1}[\rX \in N_{\bf h}({\bf x}_0)]) & = &
   {\cal V}({\bf h}) E(\XsubD \XsubD^T | \: \rX \in N_{\bf h}({\bf x}_0)) \nonumber \\
   & \equiv & {\cal V}({\bf h}) \Sigma_{1\rX}({\bf h}).
   \label{lln}
\end{eqnarray}
Because $E(\widehat \beta_1({\bf h}) \:|\: {\mathbb X}) \stackrel{P}{\longrightarrow} \beta_1({\bf h})$ 
as $n \rightarrow \infty$, this shows that the estimated local 
signal $\widehat \beta_1({\bf h})$ is conditionally
consistent given ${\mathbb X}$ as $n \rightarrow \infty$
and $d \rightarrow \infty$ when
\begin{equation}
   n^{-1} \sigma_{1\rX}^{11}({\bf h}) {\cal V}^{-1}({\bf h}) \rightarrow 0,
\end{equation}
where $\sigma_{1\rX}^{11}({\bf h})$ is defined by $(\sigma_{1\rX}^{jk}({\bf h}))
= \Sigma_{1\rX}^{-1}({\bf h}),\: 0 \leq j \leq d, \:0 \leq k \leq d$. 
%It can be shown that
%\begin{equation}
%   \sigma_{1\rX}^{11}({\bf h}) = \left({\rm Var}^{\bf h}(X_1)\left(1-\rho_1^2({\bf h})\right)\right)^{-1}
%\end{equation}
%where $\rho_1^2({\bf h})$ is the multiple correlation for the conditional 
%distribution ${\cal L}(\rX \:|\: \rX \in N_{\bf h}({\bf x}_0))$ when $X_1$ is 
%regressed on $X_2, X_3, \ldots, X_d$.
%Typically, if the variables $X_1, X_2, \ldots, X_d$ are not nearly colinear
%${\cal V}({\bf h})$ is of the form $c_1 \prod_{j=1}^d h_j$, ${\rm Var}^{\bf h}(X_1)$
%is of the form $c_2 h_1^2 \prod_{j=1}^d h_j$ and $\rho_1({\bf h}) \in (0,1)$.
%Thus, under these conditions, 
%%${\rm Var}(\widehat \beta_1({\bf h})) \rightarrow 0$ 
%conditional consistency is
%equivalent to $n^{-1} h_1^2 \rightarrow 0$ as $n,d \rightarrow \infty$.

\subsection{The Noise\label{noise}}

The denominator of the efficacy, the {\it noise}\/, 
is the asymptotic standard deviation $\sigma_0^2(\widehat \beta_1({\bf h}))$ 
of $\widehat \beta_1({\bf h})$ under $H_0$.
We derive a formula and develop a consistent estimator in what follows. 
Let $n_{\bf h}$ be the number of ${\bf X}_i$ that fall in $N_{\bf h}({\bf x}_0)$.
Note that
\begin{equation}
   \sqrt{n} \left[ \widehat \bb({\bf h}) - \bb({\bf h}) \right] = 
   \sqrt{n} \left\{ \widehat \Sigma_{1\rX}^{-1}\left[n_{\bf h}^{-1} \ldX^T {\bf Y}_{\bf h}
   - \widehat \Sigma_{1\rX} \bb\!\left({\bf h}\right) \right]\right\}
\end{equation}
where $\widehat \Sigma_{1\rX} = n_{\bf h}^{-1} (\ldX^T \ldX) \stackrel{P}{\longrightarrow}
\Sigma_{1\rX}\h$ by Equation (\ref{lln}) and
   $\left[n_{\bf h} / n {\cal V}\!\left({\bf h}\right)\right] \stackrel{P}{\longrightarrow} 1$.
By Slutsky's Theorem, $\sqrt{n}[\widehat \bb({\bf h}) - \bb({\bf h})]$ has the
same asymptotic distribution as
   $\sqrt{n} \left\{ n_{\bf h}^{-1} \Sigma_{1\rX}^{-1} \ldX^T {\bf e}_{\bf h}\right\}$
where ${\bf e}_{\bf h} \equiv {\bf Y}_{\bf h} - \bb^T\!({\bf h}) \ldX$. 
Note that by Equation (\ref{optprob}), $E^{\bf h}(\ldX^T{\bf e}_{\bf h}) = {\bf 0}$,
and with $e_i = Y_i - \bb^T\!({\bf h}){\bf X}_i$,
\begin{equation}
   n_{\bf h}^{-1} \ldX^T {\bf e}_{\bf h} =
   \left(n/n_{\bf h}\right)
   \left\{ n^{-1} \sum_{i=1}^n X_{ij} e_i {\bf 1}\!\left({\bf X}_i \in 
   N_{\bf h}\!\left({\bf x}_0\right)\right)\!: 0 \leq j \leq d\right\}.
   \label{eq3}
\end{equation}
Thus by Slutsky's Theorem and the Central Limit Theorem
\begin{equation}
   \sqrt{n} \left(n^{-1}_{\bf h} \ldX^T {\bf e}_{\bf h}\right)
   \stackrel{{\cal L}}{\longrightarrow}
   {\cal N}\!\left({\bf 0}, {\cal V}^{-1}\!\left({\bf h}\right)
   \Sigma_{{\bf Xe}}\!\left({\bf h}\right)\right)
\end{equation}
where
\begin{equation}
   \Sigma_{{\bf Xe}}\!\left({\bf h}\right) =
   \left( E\!\left(e^2 X_j X_k \:|\: {\bf X} \in N_{\bf h}\!\left({\bf x}_0\right)\right)
   \right)_{\left(d+1\right) \times \left(d + 1\right)}
\end{equation}
with $e = Y - \bb({\bf h})^T {\bf X}$.

We have shown the following.
\begin{prop}
Suppose that ${\bf h}$ and $d$ are fixed in $n$, that ${\cal V}({\bf h}) > 0$,
$0 < {\rm Var}^{\bf h}(Y) < \infty$,
and that $\Sigma_{\bf 1X}^{-1}({\bf h})$ exists; then
\begin{equation}
   \sqrt{n} \left[ \widehat \bb({\bf h}) - \bb\h\right]
   \stackrel{{\cal L}}{\longrightarrow}
   {\cal N}\!\left({\bf 0}, {\cal V}^{-1}\!\left({\bf h}\right)
   \Sigma_{1\rX}^{-1}\!\left({\bf h}\right) \Sigma_{\rX {\bf e}}\!\left({\bf h}\right)
   \Sigma_{1\rX}^{-1}\!\left({\bf h}\right)\right).
\end{equation}
%where 
%$\Sigma_{\rX {\bf e}} = {\rm Cov}(z_0^T{\bf e}, \ldots, z_d^T{\bf e}),
%\: z_j = (x_{1j}, \ldots, x_{nj})^T, \: j=0,1,\ldots, d$.
\end{prop}

%It follows that the asymptotic variance of
%$\sqrt{n}[\widehat \beta_1\h - \beta_1\h]$ 
%is
%\[
%   \sigma^2(\widehat \beta_1\h) \equiv 
%   \left[\sigma^{11}\h\right]^{-2}
%   {\rm Var}^{\bf h}(X_1 e_{\bf h})/{\cal V}\h
%\]
%where $e_{\bf h} = Y - \bb\h^T \rX$.
By using Equation (\ref{eq3}) and Liapounov's Central Limit Theorem we can
allow $d$ and ${\bf h}$ to depend on $n$ when we consider the asymptotic
distribution of $\sqrt{n}[\widehat \beta_1({\bf h}) - \beta_1({\bf h})]$.
We have shown the following.

\begin{prop}
Suppose that $\Sigma_{\bf 1X}^{-1}({\bf h})$ exists, that $[\sigma^{11}({\bf h})]^2 {\cal V}({\bf h})$
is bounded away from zero and that $0 < {\rm Var}(Y) < \infty$; then as $n \rightarrow \infty$,
\begin{equation}
   \sqrt{n}\left[\widehat \beta_1\!\left({\bf h}\right) - \beta_1\!\left({\bf h}\right)\right]
   \stackrel{{\cal L}}{\longrightarrow} {\cal N}\!\left(0, 
   \sigma^2\!\left(\widehat \beta_1\!\left({\bf h}\right)\right)\right).
\end{equation}
\label{asynorm}
\end{prop}

Because under $H_0^{(1)}$, $X_1$ is independent of $X_2, \ldots, X_d$ and $Y$,
if we set $\mu_L({\bf x}) = \beta_0 + \sum_{j=2}^d \beta_j X_j$ and 
$e = Y - \mu_L({\bf x})$,
the asymptotic variance of $\sqrt{n}[\widehat \beta_1\h - \beta_1\h]$
%conditionally on $\rX \in N_{\bf h}(x_0)$ 
is
\begin{equation}
   \sigma^2_0(\widehat \beta_1\h) \equiv
   \frac{E_{H_0}^{\bf h}\!\left[
%Y - \left(\beta_0 + \sum_{j=2}^d \beta_j X_j\right)
e
\right]^2}
   {{\cal V}\h {\rm Var}^{\bf h}_{H_0}\!\left(X_1\right)}.
\end{equation}
Now $\sigma_0(\widehat \beta_1({\bf h}))/\sqrt{n}$
is the {\it noise} part of the efficacy.
%Thus,
%\[
%   {\rm Var}_{H_0}\!\left(X_1^2 e_{\bf h}^2\right) =
%   {\rm Var}_{H_0}\!\left(X_1^2 \left[-\beta_1 X_1 + V\right]^2\right),
%\]
%where $V = Y - \sum_{j=2}^d \beta_j X_d$ is independent of $X_1$.

The sample variance $s_1^2(h_1)$ calculated using
$\{X_{i1}\!: X_{i1} \in [x_{01} - h_1, \:x_{01}+h_1]\}$
%\[
%   s_1^2\h \equiv \left(n_{\bf h} -1\right)^{-1} \sum_{\bf h}
%   \left(X_{i1} - \overline{X_1}\right)^2,
%\]
%where $\sum_{\bf h}$ stands for the sum over $i \in N_{\bf h}({\bf x}_0)$,
is our estimate of $\sigma_1^2 = {\rm Var}^{\bf h}_{H_0}(X_1)$. It is consistent
whenever $n h_1 \rightarrow \infty$ as $n \rightarrow \infty$.
%Next we consider the estimation of
%\[
%   \sigma^2_{\bf e}\h \equiv E_{H_0}^{\bf h}\!\left[Y - \mu_L\!\left(\rX\right)\right]^2
%\]
%where $\mu_L({\bf X}) = \beta_0 = \sum_{j=2}^d \beta_j X_j$ and ${\bf e} = Y-\mu_L({\bf X})$.
Note that with $\mu_0({\bf X}) = E_{H_0}^{\bf h}(Y | {\bf X})$, the null residual variance
$\sigma_e^2(\hmone) \equiv E_{H_0}^{\bf h}(e^2)$ is
\begin{equation}
   \sigma_e^2\!\left(\hmone\right)
    = \sigma^2 + E_{H_0}^{\bf h}\!\left[\mu_0\!\left(\rX\right) - \mu_L\!\left(\rX\right)\right]^2 \\
   \equiv \sigma^2 + \sigma^2_L\!\left(\hmone\right)
\end{equation}
where $\hmone=(h_2,\ldots, h_d)$ and 
$\sigma_L^2(\hmone)$ is the contribution to the variance of the error due to lack of linear fit to $\mu_0(\rX)$
under $H_0$. Let
$\widehat \mu_{L}^{\mbox{\tiny$\left(-1\right)$}}(\rX) = \widehat \beta_0 + \sum_{j=2}^d \widehat \beta_j X_j$ 
be the locally linear fit based on 
\begin{equation}
   {\cal D}\!\left(\hmone\right) \equiv
   \left\{\left(X_{ij},Y_i\right)\!: 
   X_{ij} \in \left[x_{0j}-h_j, \:x_{0j}+h_j\right],
%  X_{ij} \in N_{\bf h}({\bf x}_0), 
   \:2 \leq j \leq d,\: 1 \leq i \leq n\right\},
\end{equation}
then a natural estimator for $\sigma^2_{\bf e}$ is
\begin{equation}
   s_{\bf e}^2\!\left(\hmone\right) 
   \equiv \left(n\!\left(\hmone\right) 
   - d\right)^{-1} \sum_{{\cal D}\left(\hmone\right)} 
   \left[Y_i - \widehat \mu_L^{\supm}\!\left({\bf X}_i\right)\right]^2,
\end{equation}
where $n(\hmone)$ is the number of data points in ${\cal D}(\hmone)$.
The following can be shown (see Appendix).

\begin{lem}\label{lemma1}
If $\sigma^2 > 0$, $\sigma^2_L(\hmone) < \infty$, then, with ${\cal V}(\hmone) = E(n(\hmone))/n$,
\begin{equation}
   E_{H_0}\!\left[\sum_{{\cal D}\left(\hmone\right)} \left[Y_i - \widehat \mu_L^{\supm}\!\left({\bf X}_i\right)\right]^2 \right]
   =
   \left[ n {\cal V}\!\left(\hmone\right) -d\right] \sigma^2 + 
   n {\cal V}\!\left(\hmone\right) \sigma^2_L\!\left(\hmone\right).
\end{equation}
\end{lem}

Thus $s_e^2(\hmone)$ is a consistent estimator of $\sigma^2_e(\hmone)$ under $H_0^{(1)}$
whenever $[d/n{\cal V}(\hmone)] \rightarrow 0$ and ${\rm Var}_{H_0}[s_e^2\h] \rightarrow 0$
as $n \rightarrow \infty$ and $d \rightarrow \infty$.
It follows that under these conditions, a consistent estimate of
$\sqrt{n}\times\mbox{noise} \equiv \sigma_1({\bf h}) \equiv \sigma_0(\widehat \beta_1({\bf h}))$ is 
\begin{equation}
   \widehat \sigma_1\!\left(\hmone\right) = 
   \left\{s_{\bf e}^2\!\left(\hmone\right) \bigg/ {\cal V}\!\left(\hmone\right)
   s_1^2\!\left(h_1\right)\right\}^{1/2}.
\end{equation}
Note that because $\widehat \beta_1 \stackrel{P}{\longrightarrow} 0$ 
under $H_0^{(1)}$, we can use Slutsky's theorem to show that
this estimate is asymptotically equivalent to the usual
linear model estimate of ${\rm SD}^{\bf h}(\sqrt{n}\:\widehat \beta_1({\bf h}))$
based on data $({\bf X}_i, Y_i)$ with ${\bf X}_i \in N_{\bf h}({\bf x}_0)$.

\subsection{The Efficacy Criterion}

The usual procedures for selecting bandwidths ${\bf h}$ involve minimizing mean
squared predication or estimation error when predicting the response $Y_0$ or
estimating the mean $E(Y_0 \:|\: {\bf x}_0)$ of a case with covariate
vector ${\bf x}_0$.
The disadvantage is that significant relationships may escape procedures based
on ${\bf h}$ selected in this way because they will be based on small bandwidths
which lead to large noise.
Here we propose to choose the bandwidth to minimize the probability of Type II
error, that is, maximize power among all $t$-tests with the same asymptotic significance
level. We thereby maximize the probability of finding
significant relationships between $Y$ and a specified covariate.
This procedure automatically selects ${\bf h}$ to keep the noise small.

Because $T_n = \widehat \beta_1\h$ is asymptotically normal and $E_{H_0}(T_n) =0$,
we know (\citet*{Pitman1948}, \citet*{Pitman1979}, \citet*{Serfling1980}, \citet*{Lehmann1999})
that for contiguous alternatives the ${\bf h}$ that maximizes the asymptotic power
is the ${\bf h}$ that maximizes the absolute efficacy, where
\begin{equation}
   n^{-1/2} {\rm EFF}_1({\bf h},{\bf x}_0) = \frac{\beta_1\h}{\sigma_1\h}.
\end{equation}
%and $\sigma_1^2\h = \sigma_0^2(\widehat \beta_1\h)$. Recall $(\sigma^{ij}\h)_{d \times d}
%= \Sigma_{\rX}^{-1}\h$.
Because $\beta_1\h = \sum_{j=1}^d \sigma^{1j}\h \sigma_{Yj}({\bf h})$, 
where $\sigma_{Yj}\h = {\rm Cov}^{\bf h}(X_j, Y)$,
we can write
\begin{equation}
   n^{-1/2} {\rm EFF}_1({\bf h},{\bf x}_0) = \sum_{j=1}^d \sigma^{1j}\h \sigma_{Yj}\h / \sigma_1\h.
   \label{approxt}
\end{equation}
The efficacy optimal ${\bf h}$ for testing $H_0^{(1)}$ versus $H_1^{(1)}$ is defined by
\begin{equation}
   {\bf h}_1^{(0)} \equiv \left(h^{(0)}_{11}, \ldots, h^{(0)}_{1d}\right)^T
   = \argmax_{\bf h} {\rm EFF}_1\!\left({\bf h}, {\bf x}_0\right).
\end{equation}
For fixed alternatives that do not depend on $n$
with $\beta_1({\bf h}) >0$, this ${\bf h}_1^{(0)}$ will not satisfy
$\min\{h_{1j}^{(0)}\} \rightarrow 0$ as $n \rightarrow \infty$
because $\sigma_1({\bf h}) \rightarrow \infty$ as $\min\{h_{1j}^{(0)}\} \rightarrow 0$. 

\begin{rem}
The definition of the efficacy optimal ${\bf h}$ makes sense when the relationships
between $Y$ and covariates are monotone on $N_{\bf h}({\bf x}_0)$.
To allow for possible relationships such as
\begin{equation}
   Y = a + \cdots + \left(X_j - x_{0j}\right)^2 + \cdots + \epsilon
\end{equation}
in Section \ref{shrink} we will use neighborhoods of the form
$[x_{0j}-(1-\lambda_j)h_j, \:x_{0j}+(1+\lambda_j)h_j]$, where $-1 \leq \lambda_j \leq 1$,
rather than the symmetric intervals $[x_{0j}-h_j, \:x_{0j} + h_j]$. 
%See Section \ref{???}.
The properties of the procedures are not changed much by the introduction of the
extra tuning parameters $\lambda_1, \lambda_2, \ldots, \lambda_d$. See \citet*{Doksum2006}
for the case $d=1$.
\end{rem}

The estimate of $n^{-1/2} {\rm EFF}_1$ is $n^{-1/2} t_1({\bf h}, {\bf x}_0)$ where
$t_1({\bf h}, {\bf x}_0)$ is the $t$-statistic
\begin{equation}
   t_1\!\left({\bf h},{\bf x}_0\right) = \left\{
   \begin{array}{ll}
%   \frac{{\cal V}\!\left({\bf h}\right) \widehat \beta_1\!\left({\bf h}\right)
%   s_1\!\left({\bf h}\right)}{s_{\bf e}\!\left({\bf h}\right)} &
    \frac{n^{1/2}\!\left(\hmone\right) \widehat \beta_1\!\left(\hmone\right)
    s_1\!\left(h_1\right)}{n^{1/2} s_{\bf e}\!\left({\bf h}\right)} &
   \mbox{if $n\!\left(\hmone\right) > d+1$} \\
   0 & \mbox{otherwise}.
   \end{array}
   \right.
\end{equation}
Write $d_n$ for $d$ to indicate that $d$ may depend on $n$.
Using the results of Sections \ref{signal} and \ref{noise}, we have the following.

\begin{prop}\label{consest}
In addition to the assumptions of Proposition \ref{asynorm}, assume that
as $n \rightarrow \infty$,
$[d_n/n {\cal V}({\bf h})] \rightarrow 0$, 
%$n^{-1} \sigma_{1x}({\bf h}) {\cal V}^{-1}({\bf h}) \rightarrow 0$,
$[\sigma^{11}({\bf h})]^2 {\cal V}({\bf h})$ is bounded away from zero,
and ${\rm Var}_{H_0}[s^2_{\bf e}({\bf h})] \rightarrow 0$; then $t_1({\bf h},{\bf x}_0)$ is a consistent
estimate of ${\rm EFF}_1$, in the sense that
\begin{equation}
   \frac{t_1\!\left({\bf h},{\bf x}_0\right)}{{\rm EFF}_1\!\left({\bf h},{\bf x}_0\right)}
   \stackrel{P}{\longrightarrow} 1 \:\:\: \mbox{as $n \rightarrow \infty$}.
\end{equation}
\end{prop}

\begin{coro}\label{coro1}
Suppose ${\cal G}$ (a grid) denotes a finite set of vectors of the form $({\bf h},{\bf x}_0)$,
then, under the assumptions of Proposition \ref{consest},
 the maximizer of $t_1({\bf h},{\bf x}_0)$ over ${\cal G}$ is asymptotically optimal in
the sense that
\begin{equation}
   \frac{
   \max \left\{t_1\!\left({\bf h},{\bf x}_0\right)\!: \left({\bf h},{\bf x}_0\right) \in {\cal G}\right\}}
   {
   \max \left\{{\rm EFF}_1\!\left({\bf h},{\bf x}_0\right)\!: \left({\bf h},{\bf x}_0\right) \in {\cal G}\right\}}
   \stackrel{P}{\longrightarrow} 1.
\end{equation}
\end{coro}

%\begin{rem}
%One grid search for ${\bf h}$ that works well is the following: For each $j$, standardize $X_{ij},
%\:1 \leq i \leq n$ to have SD equal to one. Then set ${\bf h}^{(1)} = (h^{(1)}, \ldots, h^{(1)})$
%where $h^{(1)}$ is selected so that $N_{{\bf h}^{(1)}}({\bf x}_0)$ contains 25\% of the observations.
%Next maximize $t_j({\bf h}^{(1)}, {\bf x}_0)$ with respect to the $j^{th}$ coordinate keeping
%the other coordinates equal to $h^{(1)}$. This yields a new ${\bf h}^{(2)} = (h_1^{(2)}, \ldots, h_d^{(2)})$.
%Now maximize $t_j({\bf h}^{(2)}, {\bf x}_0)$ with respect to the $j^{th}$ coordinate keeping the
%other coordinates fixed at $h_k^{(2)}, \: k \neq j$. This yields ${\bf h}^{(3)} = (h_1^{(3)}, \ldots, h_d^{(3)})$;
%then stop.
%\end{rem}

\begin{rem}
\citet*{Hall2000} used the maximum of local $t$-statistics to test the global hypothesis that $\mu(\cdot)$ is monotone
in the $d=1$ case. Their estimated local regression slope is the least squares estimate based on $k$ nearest neighbors
and the maximum is over all intervals with $k$ at least 2 and at most $m$. They established unbiasedness and consistency
of their test rule under certain conditions.
\end{rem}

Because the asymptotic power of the test based on $\widehat \beta_1({\bf h})$ tends to one for all
${\bf h}$ with $\beta_1({\bf h}) >0$ and $|{\bf h}| > 0$, it is more interesting to consider
Pitman contiguous alternatives $H_{1n}^{(1)}$ where $\beta_1({\bf h})$ depends on $n$ and
tends to zero at a rate  that ensures that the limiting power is between $\alpha$ and 1, where
$\alpha$ is the significance level. That is, we limit the parameter set to the set where deciding
between $H_0^{(1)}$ and the alternative is difficult. We leave out the cases where the right decision
will be reached for large $n$ regardless of the choice of ${\bf h}$.

Now the question becomes: For sequences of contiguous alternatives with $\beta_1({\bf h}) = \beta_{1n}({\bf h})
\rightarrow 0$ as $n \rightarrow \infty$, what are the properties of $h_{1n}^{(0)}$? In
particular, does the efficacy optimal $h_{1n}^{(0)}$ tend to zero
as $n \rightarrow \infty$? The answer depends on the alternative $H_{1n}^{(1)}$ as will be
shown below.

\subsection{Optimal Bandwidths for Pitman Alternatives with Fixed Support}
Under $H_0^{(1)}$,
\begin{equation}
   \beta_1\!\left({\bf h}\right) = \sum_{j=1}^d \beta_{1j}\!\left({\bf h}\right)
   \equiv \sum_{j=1}^d \sigma^{1j}\!\left({\bf h}\right) {\rm Cov}^{\bf h}\!\left(X_j,Y\right)
   = 0.
\end{equation}
We consider sequences of contiguous Pitman alternatives with $\beta_{1j}({\bf h}) \propto cn^{-1/2}$,
such as
\begin{equation}
   Y = \alpha + \gamma_n r\!\left({\bf X}\right) + \epsilon
   \label{pitalt}
\end{equation}
where $\gamma_n \equiv c n^{-1/2}$, $c \neq 0$, and $|{\rm Cov}^{\bf h}(X_j, r({\bf X}))| > b$,
for $b>0$. 
%is of this form provided $r({\bf x})$ and its support $\{{\bf x}\!: |r({\bf x})| >0\}$
%does not depend on $n$. For such alternatives
Here,
\begin{equation}
   {\rm EFF}_1\!\left({\bf h},{\bf x}_0\right)
   \longrightarrow
   \left(
   \frac{c {\cal V}\!\left({\bf h}\right) {\rm SD}^{h_1}\!\left(X_1\right)}
%  {\sigma^2 + \sigma_L^2\!\left({\bf h}\right)}\right)
   {\sigma_e^2\!\left(\hmone\right)}\right)
   \sum_{j=1}^d \beta_{1j}\!\left({\bf h}\right)
\end{equation}
with $\beta_{1j}({\bf h}) = \sigma^{1j}({\bf h}) {\rm Cov}^{\bf h}(X_j,r({\bf X}))$.
As in the case of fixed alternatives,
the maximizer ${\bf h}_1^{(0)}$ does not satisfy $\min\{h_{1j}^{(0)}\} \rightarrow 0$
as $n \rightarrow \infty$ because
${\rm EFF}_1({\bf h},{\bf x}_0) \rightarrow 0$ as $\min\{h_{1j}^{(0)}\} \rightarrow 0$.
This ${\bf h}_1^{(0)}$ is asymptotically optimal in the sense of Corollary \ref{coro1}.

\subsection{Comparison with Mean Squared Error}

There is a large literature on selecting bandwidths by minimizing mean squared
error (MSE). Here MSE can be expressed as the following.

\begin{prop}
If $0 < {\rm Var}(Y) < \infty$, then
\begin{eqnarray}
   E^{{\bf h}}\!\left\{ \left[
   \widehat \mu_L\!\left({\bf x}_0\right) - \mu\!\left({\bf x}_0\right)\right]^2
   |\: {\mathbb X}
   \right\}
   & = & \frac{\sigma^2}{n_{\bf h}} + \left[\mu_L\!\left({\bf x}_0\right) - \mu\!\left({\bf x}_0\right)
   \right]^2 \nonumber \\
   & = &
   \frac{\sigma^2}{n {\cal V}\!\left({\bf h}\right)}
   + \left[ \mu_L\!\left({\bf x}_0\right) - \mu\!\left({\bf x}_0\right)\right]^2 
   + o_P\!\left(1/n\right).
   \label{approxmse}
\end{eqnarray}
\end{prop}

\begin{proof}
MSE is variance plus squared bias where the conditional squared bias is as given.
The conditional variance of the local least squares estimate $\widehat \mu_L({\bf x}_0) = \widehat \beta_0({\bf h})$
given ${\mathbb X}$ is $\sigma^2/n_{\bf h}$ where $(n_{\bf h}/n) \stackrel{P}{\rightarrow} {\cal V}({\bf h})$.
\end{proof}

%We next compare (\ref{approxt}) with the mean squared error,

%\begin{prop}
%\noindent {\em Proposition 1.2.} 
%Under the preceding assumptions and those of Appendix I,
%\begin{equation}
%   E\!\left[
%   \left(\widehat \mu_L\!\left({\bf X}\right) - 
%   \mu\!\left({\bf X}\right)\right)^2
%   1_{\{{\bf X} \in N_{\bf h}\!\left({\bf x}_0\right)\}}
%   \right]
%   = \frac{\left(d+1\right) \sigma^2}{n{\cal V}\!\left({\bf h}\right)} + \sigma_L^2\!\left({\bf h}\right)
%   + o_P\!\left(1\right).
%   \label{approxmse}
%\end{equation}
%\end{prop}

If finding significant dependencies is the goal, we
prefer maximizing Equation (\ref{approxt}) to minimizing Equation (\ref{approxmse})
because
the local MSE (\ref{approxmse}), as well as its global version, focuses on finding the
${\bf h}$ that makes $\widehat \mu_L({\bf x})$ close to the true unknown curve $\mu({\bf x})$.
Using results of \citet*{Ruppert1994} and \citet*{Fan1996}, page 302,
we can show that under regularity conditions (the Hessian of $\mu({\bf x})$ exists), the bandwidths
minimizing (\ref{approxmse}) tend to zero at the rate $n^{-1/(d+4)}$. By plugging such bandwidths
into Equation (\ref{approxt})
we see
that this will in many cases,
make it nearly impossible to find dependencies. Various semiparametric assumptions
have brought the rate of convergence to zero of the bandwidths to $n^{-1/5}$ which still makes it
likely that dependencies will be missed. By using (\ref{approxt}) we have a simple method for finding
bandwidths that focus on finding dependencies for any type of alternative. However there are 
alternatives where ${\bf h} \rightarrow {\bf 0}$ makes sense. We construct these next.

\section{Bandwidths for Alternatives with Shrinking Support}\label{shrink}

We consider sequences of alternatives where the set $A_n =\{{\bf x}\!: |\mu({\bf x})-\mu_Y| > 0\}$ is
a connected region whose volume tends to zero as $n \rightarrow \infty$.
% BEGIN INSERT FROM OTHER PAPER
One such set of alternatives is given by Pitman alternatives of the form
\begin{equation}
   K_n\!:\:\: Y = a + \sum_{j=1}^d \gamma_j \:W_j\!\left(\frac{X_j-x_{0j}}{\theta_j}\right) +\epsilon,
   \label{Pitmod}
\end{equation}
where ${\bf X}$ and $\epsilon$ are uncorrelated, $\epsilon$ has mean zero and variance $\sigma^2$,
and each $W_j(\cdot)$ has support $[-1,1]$. 
We assume that each $X_j$ is continuous, in which case
the hypothesis holds
with probability one when $\gamma_j \theta_j = 0$ for all $j$.
We consider models where $\theta_j = \theta_j^{(n)} \rightarrow 0$ as
$n \rightarrow \infty$, and $\gamma_j$ may or may not depend on $\theta_j$ and $n$. 
For these alternatives
the neighborhood where $E(Y|{\bf X})$ is non-constant shrinks to achieve a Pitman balanced
model where the power converges to a limit between the level $\alpha$ and 1 as
$n \rightarrow \infty$. Note, however, that the alternative does not depend on ${\bf h}$.
We are in a situation where ``nature'' picks the neighborhood size ${\bf \theta} \equiv (\theta_1, \ldots,\theta_d)^T$, 
and
the statistician picks the bandwidth ${\bf h}$. This is in contrast to 
\citet*{Blyth1993}, \citet*{Fan2001}, and \citet*{Ait2001}
%Blyth (1993) and A\"{\i}t-Sahalia, Bickel and Stoker (2001) 
who let the alternative depend on ${\bf h}$.

We next show that for fixed ${\bf h}$ with $|{\bf h}|>0$, ${\rm EFF}_1({\bf h},{\bf x}_0) \rightarrow 0$
as $\max_j\{\gamma_j\theta_j\} \rightarrow 0$ in model (\ref{Pitmod}).
First note that
\begin{equation}
   {\rm Cov}^{\bf h}\!\left(X_j, Y\right) = \sum_{k=1}^d \gamma_k {\rm Cov}^{\bf h}\!\left(
   X_j,W_k\!\left(\frac{X_k-x_{0k}}{\theta_k}\right)\right).
\end{equation}
Next note for $h_k >0$ fixed, $\gamma_k$ bounded from above, and $\theta_k \rightarrow 0$,
\begin{equation}
   \gamma_k {\rm Cov}^{\bf h}\!\left(X_j, W_k\!\left(\frac{X_k - x_{0k}}{\theta_k}\right)\right)
   \rightarrow 0
\end{equation}
because $W_k((X_k-x_{0k})/\theta_k)$ tends to zero in probability as $\theta_k \rightarrow 0$ and
any random variable is uncorrelated with a constant. This heuristic can be verified (see the appendix)
by a change of variable. This result is stated more precisely as follows.

\begin{prop}\label{prop5}
If ${\bf h}$ is fixed with $|{\bf h}|>0$, then as $\max_k \{\gamma_k\theta_k\}\rightarrow 0$
in model (\ref{Pitmod}), \\
(a) ${\rm EFF}_1\!\left({\bf h},{\bf x}_0\right) = O\!\left(\max_k\{\gamma_k \theta_k\}\right)$ and \\
(b) $t_1\!\left({\bf h},{\bf x}_0\right) = O_P\!\left(\max_k\{\gamma_k \theta_k\}\right)$.
\end{prop}

Proposition \ref{prop5} shows that for model (\ref{Pitmod}), fixed ${\bf h}$ leads
to small ${\rm EFF}_1({\bf h},{\bf x}_0)$.
Thus we turn to the $h_1 \rightarrow 0$ case. If $h_1 > 0$, then observations $(X_1,Y)$ with $X_1$ outside
$[x_{01}-\theta_1, x_{01}+\theta_1]$ do not contribute to the estimation of $\beta_1({\bf h})$. Thus choosing
a smaller $h_1$ may be better even though a smaller $h_1$ leads to a larger variance for $\beta_1({\bf h})$.
This heuristic is made precise by the next results which provide conditions where $h_1=\theta_1$ is the
optimal choice among $h_1$ satisfying $h_1 \geq \theta_1$.
First, define
%\begin{equation}
%   {\cal V}_{-1} \equiv P\!\left(X_j \in \left[ x_{0j}-h_j, x_{0j}+h_j\right], \: 2 \leq j \leq d\right)
%\end{equation}
%and
\begin{equation}
   m_j\!\left(W_1\right) \equiv \int_{-1}^1 s^j W_1\!\left(s\right) ds, \:\:\:\:\:\: j=0, 1, 2.
\end{equation}

\begin{thm}\label{thm1}
Assume that $X_1$ is independent of $X_2,\ldots,X_d$, that the density $f_1(\cdot)$ of $X_1$
has a bounded, continuous second derivative at $x_0$, and that $f_1(x_0) >0$. Then, in model
(\ref{Pitmod}), as $\theta_1 \rightarrow 0$, $h_1 \rightarrow 0$ with $h_1 \geq \theta_1$,
and $\hmone$ fixed with $|\hmone|>0$,\\
(a) ${\rm Cov}^{\bf h}\!\left(X_1,Y\right) = \gamma_1 \theta_1^2 m_1\!\left(W_1\right)\bigg/ 2h_1
   + o\!\left(\frac{\gamma_1 \theta_1^2}{h_1}\right)$.\\
(b) If $m_1(W_1) \neq 0$ and $m_2(W_1) \neq 0$, then
\begin{equation}
   n^{-1/2} {\rm EFF}_1\!\left({\bf h}, {\bf x}_0\right) 
   \propto
   \left[\sigma^2 + \sigma_{\bf L}^2\!\left(\hmone\right)\right]^{1/2}
   \gamma_1 \theta_1^2 h_1^{-3/2} m_1\!\left(W_1\right) f_1^{1/2}\!\left(x_{01}\right)
%  \prod_{j=2}^d P\!\left(x_{0j}-h_j \leq X_j \leq x_{0j}+h_j\right)
   {\cal V}^{1/2}\!\left(\hmone\right)
\end{equation}
which is maximized 
%by taking $h_j = \infty$, $j \geq 2$, and maximized 
subject to $h_1 \geq \theta_1$ by $h_1 = \theta_1$.
\end{thm}

\begin{thm}\label{thm2}
Suppose the joint density $f({\bf x})$ of ${\bf X}$ has a bounded continuous Hessian
at ${\bf x}_0$, and that $f({\bf x}_0) >0$, then in model (\ref{Pitmod}), as $|{\bf h}| \rightarrow 0$
and $|\btheta| \rightarrow 0$ with $h_1 \geq \theta_1$, conclusion (b) of Theorem \ref{thm1}
with ${\cal V}^{1/2}(\hmone)$ replaced by $\{2^{d-1} \prod_{j=2}^d h_j\}^{1/2}$ holds.
\end{thm}

The case where $h_1 < \theta_1$ remains; then some of the data in the neighborhood where
the alternative holds is ignored and efficacy is reduced. Under suitable conditions on $W_j$, $j=1,\ldots,d$, the
optimal $h_1$ equals $\theta_1$; that is, $2h_1^{(0)}$ is the length of the interval where the $X_1$ term in the additive
model (\ref{Pitmod}) is different from zero. 
Details of the justification can be constructed by extending the results of \citet*{Doksum2006} to the $d>1$ case.

\begin{rem}
\citet*{Fan1992}, \citet*{Fan1993}, and \citet*{Fan1996}
considered minimax kernel estimates for models where $Y= \mu(X) + \epsilon$ in the $d=1$ case
and found that the least favorable distribution for estimation of $\mu(x_0)$ using asymptotic 
mean squared error has $Y=\mu_0(X) +\epsilon$ with
\begin{equation}
   \mu_0\!\left(x\right) = \frac{1}{2} b_n^2 \!\left[
   1 - c \left(\frac{x-x_0}{b_n}\right)^2
   \right]_+
\end{equation}
where $b_n = c_0 n^{-1/5}$, for some positive constants $c$ and $c_0$. Similar results were 
obtained for the white noise model by \citet*{Donoho1991a,Donoho1991b}. 
\citet*{Lepski1999}, building on \citet*{Ingster1982}, considered minimax
testing using kernel estimates of $\mu(x)$ and a model with normal
errors $\epsilon$ and ${\rm Var}(\epsilon) \rightarrow 0$ as $n \rightarrow \infty$.
Their least favorable distributions (page 345) has $\mu_0(x)$ equal to a random linear combinations of
functions of the form
\begin{equation}
   \mu_j\!\left(x\right) = \left(h^{1/2} \!\int W^2\!\left(t\right) dt\right)^{-1}
   W\!\left(\frac{x-t_j}{h}\right).
\end{equation}
%which seems to indicate a model that depends on $h$. Here each $\mu_j(x)$ has
%$\|\mu_j\|^{-\infty} \rightarrow 0$ as $n \rightarrow \infty$.
\end{rem}

\section{Variable Selection}\label{varsel}

Suppose we want to test $H_0^{(1)}$ that $X_1$ and $Y$ are unrelated and wonder whether
we should keep $X_j$, $j \geq 2$, in the model when we construct a $t$-statistic for this testing
problem. In experiments where confounding variables are possibly present, 
it would seem unreasonable to keep or exclude $X_j$ on the basis of power
for testing $H_0^{(1)}$ because confounding could lead to dropping $X_j$ in situations where
the relationship between $X_1$ and $Y$ is spurious. For this reason it would seem more
reasonable to base the decision about keeping $X_j$ on the strength of its relationship
to $(X_1,Y)$ and thus we ask whether $X_j$ contributes to accurate
prediction of the response $Y$ conditionally given $X_1$.
More generally we consider dropping several variables and consider fits
to $Y$ over subsets of $\Re^d$ that vary both in dimension and size (tuning parameter) within
each dimension and we select the dimension and size which maximizes efficacy in the $X_1$
direction while controlling for spurious correlation.

We define a {\it feature} ${\cal S}_k$ 
to be a subset of the covariate space which yields the maximal
absolute $t$-statistic for some ${\bf x} \in {\cal S}_k$. In other words, among all of the
subregions over which linear models are fit, we discard those that do not
maximize the absolute $t$-statistic for {\bf some} covariate vector ${\bf x}$. The
remaining regions are denoted ${\cal S}_1, {\cal S}_2, \ldots, {\cal S}_r$.
Assume the features are placed in decreasing order
of the value of absolute $t$-statistic for the variable of interest.
These features are subsets of $\Re^d$ ordered with respect to their relevance to the relationship
between $X_1$ and $Y$.
Define
\begin{equation}
   {\cal S}_k' \equiv {\cal S}_k \cap \left({\cal S}_1' \cup {\cal S}_2' \cup \cdots \cup {\cal S}_{k-1}'\right)^c
\end{equation}
with ${\cal S}_1' \equiv {\cal S}_1$. Then $\cup_k {\cal S}_k' = \cup_k {\cal S}_k$ and the
${\cal S}_k'$ are disjoint.

%For illustrative purposes, 
%assume that there are two covariates, $X_1$ and $X_2$.
%Define
%\begin{equation}
%   \mu_i'\!\left(x_1, x_2\right) \equiv 
%   \left.\frac{d\mu\!\left(x'_1,x'_2\right)}{d x'_i}\right|_{(x_1,\:x_2)}
%\end{equation}
%for $i=1,2$. Assume that we are interested in the relationship
%between $X_1$ and $Y$ when $X_1 \approx x_1$. As $x_2$ is allowed to
%vary, $(x_1,x_2)$ will either remain in the same subset ${\cal S}_k'$,
%or move between multiple such subsets. The latter will be the case
%if either (1) $\mu_1'(x_1,x_2)$ varies with $x_2$ or (2) $\mu_2'(x_1,x_2)$
%varies with $x_2$, or both. If only case (1) were true, then the
%covariate $X_2$ should be removed as a term in the local models fit around $X_1 = x_1$.
%But, if case (2) were true, removing $X_2$ could lead to spurious correlation
%between $X_1$ and $Y$, if $X_1$ and $X_2$ are related.
%Apparent correlation between $X_1$ and $Y$ is spurious if that correlation
%disappears once the model is fit with additional covariates.

There are competing goals: We want the model to be parsimonious (not
include too many covariates), especially since we are focusing on cases
where $d$ is large. But, we don't want to exclude any
covariate which, if included, changes the picture of the $(X_1,Y)$
relationship.
Thus if $X_1$ and $X_2$ are closely related,
as are $X_2$ and $Y$, there will also be an apparent relationship
between $X_1$ and $Y$. 
%If the goal were simply to construct a model
%which leads to accurate predictions of $Y$, then it might be sufficient to
%include only $X_1$, since $X_2$ may not add any predictive power.
%Here, we take a more conservative approach, and imagine that our primary
%goal is 
To avoid making false claims regarding the strength of the relationship
between $X_1$ and $Y$, we tentatively include $X_2$, and let the analysis
show whether that weakens the $(X_1,Y)$ relationship.

This analysis of the $(X_1,Y)$ relationship after correcting
for other $X's$
can be accomplished in the following way: Consider the quantity
\begin{equation}
   \gamma(x_1,\widehat \mu_j^{({\cal S})}) \equiv E\!\left[\left(Y - 
   \widehat \mu_j^{({\cal S})}({\bf X})\right)^2 \:|\: X_1 = x_1\right],
\end{equation}
the expected squared prediction error when $X_1 =x_1$ , where $({\bf X},Y)$
is independent of $({\bf X}_i,Y)$, $1\leq i\leq n$, and
$\widehat \mu_j^{({\cal S})}({\bf x})$ stands for the model fit for
model $j$ based on a subset of covariates that are restricted to
subregion ${\cal S}$.
We seek
the subset of the covariates which minimizes this criterion. 
If there is no relationship
between $X_2$ and $Y$ when $X_1 = x_1$,
then $X_2$ will be excluded based on comparisons of $\gamma(x_1,\widehat \mu_j^{({\cal S})})$
for different models $j$;
this could lead to different covariates included for different values of $X_1$.

%Just as we asked questions regarding the power against $H_0^{(1)}$ for different values of
%$X_1$, and then chose the neighborhoods in an effort to maximize that power, here
%we look at the mean squared prediction error as a function of $X_1$.
%This leads us to consider the quantity
%\[
%   \gamma(x,\widehat \mu^{({\cal S})}) \equiv E\!\left[\left(Y - 
%   \widehat \mu^{({\cal S})}({\bf X})\right)^2 \:|\: X_1 = x\right],
%\]
%the expected prediction error when $X_1 = x$, where $({\bf X},Y)$ is
%independent of the sample $({\bf X}_i, Y_i), 1 \leq i \leq n$.
%For example, imagine
%two competing models, with corresponding model fits
%$\widehat \mu_1$ and $\widehat \mu_2$, and suppose these models
%differ in that the former includes covariate $X_j$, and the latter does not.
%If $\gamma(x,\widehat \mu_1) < \gamma(x,\widehat \mu_2)$, then 
%it is important to include $X_j$ in the local model, since if it were not included there would be
%a large amount of variability in $Y$ that would be left uncorrected for, and our picture of
%the relationship between $X_1$ and $Y$ would be blurred.
%Thus, the variables selected for inclusion in the
%model can be different for different values of $X_1$ since the importance of keeping $X_j$ 
%in the model for predicting $Y$ may depend on $X_1$. 

This is implemented as follows.
Within each ${\cal S}_k$, $m$ different linear models are fit; these models
differ in which covariates are included, but the notation is consistent in the sense that ``model $j$'' 
always refers to the same list of covariates. 
Let $\widehat \mu_{jk}({\bf x})$ 
denote the fitted value at covariate value ${\bf x} \in {\cal S}_k$ from model $j$ fit over ${\cal S}_k$ and let
$\widehat \mu_j({\bf x}) \equiv \widehat \mu_{jk}({\bf x})$ if ${\bf x} \in {\cal S}_k'$. Each ${\bf x}$ lies in exactly
one ${\cal S}_k'$, so this is uniquely defined.

The mean squared prediction error for a model fit
%\begin{equation}
%  R\!\left(\widehat \mu\right) \equiv E\!\left[\left(Y - \widehat \mu({\bf X})\right)^2\right] 
%\end{equation}
is commonly approximated by the leave-one-out cross validation score
%\begin{equation}
%   \widehat R\!\left(\widehat \mu\right) \equiv 
%   n^{-1} \sum_{i=1}^n \left(Y_i - 
%   \widehat \mu_{-i}({\bf X}_i)\right)^2
%\end{equation}
based on the fitted value
$\widehat \mu_{-i}({\bf X}_i)$ from the fit with $({\bf X}_i,Y_i)$ removed. If
$\widehat \mu$ is a linear model fit using least squares then the cross-validation
prediction error is
\begin{equation}
   \left(Y_i - 
   \widehat \mu_{-i}({\bf X}_i)\right)^2
   =
   \left(\frac{Y_i - \widehat \mu({\bf X}_i)}{1-h_i} \right)^2
\end{equation}
where $h_i$ is the $i^{th}$ diagonal element of the hat matrix from the linear model.

Here, we do not fit a global linear model, but the fitted value $\widehat \mu_j^{({\cal S})}({\bf X}_i)$ 
is the result of the fit from {\bf some} linear model. Let $h_{ij}$ denote the diagonal element corresponding
to ${\bf X}_i$ of the hat matrix from the linear model fit which 
gives us $\widehat \mu_j^{({\cal S})}({\bf X}_i)$. Then define
\begin{equation}
   \widehat \gamma_{ij}^{({\cal S})} 
   \equiv \left(\frac{Y_i - \widehat \mu_j^{({\cal S})}({\bf X}_i)}{1-h_{ij}} \right)^2.
\end{equation}
%so that
%\begin{equation}
%   \widehat R\!\left(\widehat \mu_j\right) = n^{-1} \sum_{i=1}^n \gamma_{ij}.
%\end{equation}
Note how this automatically adjusts for the differing degrees of freedom in the
different models since $\sum (1-h_{ij}) = n-d-1$.
We can now estimate $\gamma(\cdot, \widehat \mu_j^{({\cal S})})$ 
by the smooth $\widehat \gamma(\cdot,\widehat \mu)$ of $\widehat \gamma_{ij}^{({\cal S})}$ versus
the observed value of $X_1$ for each data point ${\bf X}_i, i=1,2,\ldots, n$. 
For a given $X_1 = x_1$, the model $j$ within ${\cal S}_k$ 
with the smallest $\widehat \gamma(x_1, \widehat \mu_j^{({\cal S}_k)})$ is selected.
Figure \ref{case1_cvplot} shows an example where four models are compared.

Note that we avoid the curse of dimensionality by using features ${\cal S}_k$ with large values
of the $t$-statistics in direction $X_1$. If there are a large number of variables, instead of
considering all possible combinations of variables, we can use the backward deletion approach
commonly used in connection with AIC or SBC (Schwarz's Bayesian Criteria).

To summarize, by selecting the features as described above, we are simultaneously selecting the
number of variables to include and the size of candidate neighborhoods for computing $t$-statistics
in a given direction, here $X_1$. We select the neighborhoods where the $t$-statistics are maximized
and we use conditional prediction error to select the variables
in such a way that we are protected against using models that produce spurious correlations.

\begin{rem}
A great number of tests of model assumptions and variable selection procedures
are available in a nonparametric setting, e.g.
\citet*{Azzalini1989},
\citet*{Raz1990},
\citet*{Eubank1993},
\citet*{Hardle1993},
\citet*{Bowman1996},
\citet*{Hart1997},
\citet*{Stute1997},
\citet*{Lepski1999},
\citet*{Fan2001},
\citet*{Polzehl2002},
\citet*{Zhang2003b},
\citet*{Samarov2005},
among others.
One class of tests
of whether the $j$th variable has an effect
looks at the difference between the mean $\mu ({\bf x})$ of $Y$ given all the ${x}$'s 
and the mean $\mu_{-j} ({\bf x}_{-j})$ of $Y$ given all but the $j$th variable.  
Specifically, for some weight functions $w({\bf X})$, measures of the form
\[   
   m_j = E \{ [ \mu ({\bf X}) - \mu_{-j} ({\bf X}_{-j}) ]^2 w({\bf X}) \}
\]
are considered (\citet*{Doksum1995}, \citet*{Ait2001}).  
Similar measures compare prediction errors of the form
\[
   E\{ [ Y - \mu ({\bf X}) ]^2   w({\bf X})\} \:\:\:\mbox{and}\:\:\:
   E\{ [ Y - \mu_{-j} ({\bf X}_{-j} )]^2 w({\bf X}_{-j} ) \}.
\]
Our efficacy measure  is a version of $\mu({\bf x}) - \mu_{-j}({\bf x}_{-j})$ adjusted
for its estimability while our measure of spurious correlation is a conditional version
of the above prediction errors.
\end{rem}

\section{Critical Values}\label{critvals}

\begin{figure}[ht]
\begin{center}
\psfig{file=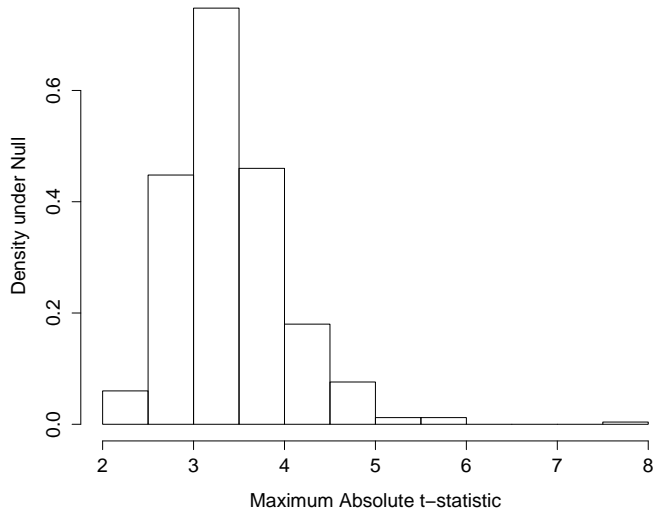,width=\onecolwidth}
\caption{The simulated distribution of the maximum absolute $t$-statistic at ${\bf x}_0 = (0.4, 0.4,0.4)$
using the model described in Section \ref{mc}, Equation (\ref{simpmod}).}
\label{permdist}
\end{center}
\end{figure}

First consider testing $H_0^{(1)}\!:$``$X_1$ is
independent of $X_2, X_3, \ldots, X_d, Y$'' against
the alternative $H_1^{(1)} = H_1^{(1)}({\bf x}_0)$ that
$\mu({\bf x})$ is not constant as a function of $x_1$ in
a neighborhood of a given covariate vector ${\bf x}_0$
from a targeted unit of interest. Our procedure uses the
data to select the most efficient $t$-statistic adjusted
for spurious correlation, say $t_1({\bf x}_0, ({\bf X}, {\bf Y}))$,
where $({\bf X},{\bf Y}) = ((X_{ij})_{n\times d}, (Y_i)_{n\times 1})$
are the data. Let ${\bf X}_1^*$ be a random permutation of $(X_{11},\ldots, X_{n1})$
and ${\bf X}^* = ({\bf X}_1^*, {\bf X}_2, \ldots, {\bf X}_d)$ where 
${\bf X}_j = (X_{1j},\ldots, X_{nj})^T$. Then, for all $c>0$,
\[
   P_{H_0}\!\left(|t_1\!\left({\bf x}_0, \left({\bf X},{\bf Y}\right)| \leq c\right)\right)
   = P\!\left(|t_1\!\left({\bf x}_0, \left({\bf X}^*,{\bf Y}\right)| \leq c\right)\right).
\]
Next, we select $B$ independent random permutations ${\bf X}_1^*, \ldots, {\bf X}_B^*$
and use the $(1-\alpha)$ sample quantile of $|t_1({\bf x}_0, ({\bf X}_k^*,{\bf Y}))|$,
for $k=1,2,\ldots, B$, as the critical value.
As $B \rightarrow \infty$, this quantile converges in probability to the level $\alpha$
critical value.
Figure \ref{permdist} gives as example of 
the simulated distribution of $|t_1({\bf x}_0, ({\bf X}, {\bf Y}))|$ under
the null hypothesis for ${\bf x}_0 = (0.4, 0.4, 0.4)$.
Here $B=500$, and the model is described in Section \ref{mc}, stated in
Equation (\ref{simpmod}). 

Note that although $t_1({\bf x}_0, ({\bf X},{\bf Y}))$ is selected using absolute values
of efficacy, we can still perform valid one-sided tests by following the above procedures
without the absolute values. In this case the alternative is that $\mu({\bf x})$ is increasing
(or decreasing) as a function of ${\bf x}_1$ in a neighborhood of ${\bf x}_0$.

Next consider testing $H_0^{(1)}$ against the alternative $H_{(1)}$ that $X_1$ is not
independent of $X_2,\ldots, X_d, Y$. If ${\bf x}^{(1)}, \ldots, {\bf x}^{(g)}$ is a set
of grid points we can use the sum, sum of squares, maximum, or some other norm, of
\[
   \left|t_1\!\left({\bf x}^{(1)}, \left({\bf X},Y\right)\right)\right|, \ldots,
   \left|t_1\!\left({\bf x}^{(g)}, \left({\bf X},Y\right)\right)\right|
\]
to test $H_0^{(1)}$ versus $H_{(1)}$. We may also instead of grid points use ${\bf X}_1,
\ldots, {\bf X}_n$ or a subsample thereof. Again the permutation distribution of these
test statistics will provide critical values.

Finally, consider the alternative that $\mu({\bf x})$ is monotone in $x_1$ for all $x_1 \in \Real$.
We would proceed as above without the absolute values, use a left sided test for monotone decreasing,
and a right sided test for monotone increasing, see \citet*{Hall2000} who considered the $d=1$ case.

\section{The Analysis Pipeline; Examples}\label{mc}

\begin{figure}[ht]
\begin{center}
\psfig{file=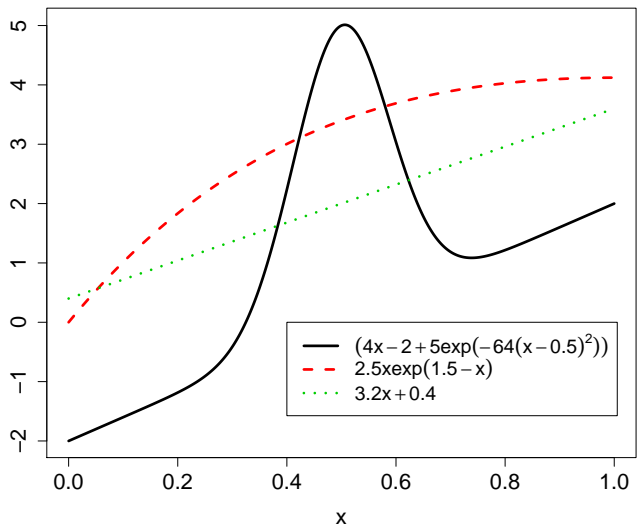,width=\onecolwidth}
\caption{The three functions which will be used in the simulations.}
\label{curveplot}
\end{center}
\end{figure}

In this section we describe practical implementation of the data analysis 
pipeline developed in the previous sections, using both simulated and real
data.
%See Figure \ref{pipeline}.
The simulations are built around the following functions, plotted in Figure
\ref{curveplot}.
\begin{eqnarray*}
f_1(x) & = & 4x-2 + 5\exp(-64(x-0.5)^2) \\
f_2(x) & = & 2.5x \exp(1.5-x) \\
f_3(x) & = & 3.2x + 0.4 
\end{eqnarray*}

The algorithm can be described as follows.
Fix a covariate vector ${\bf x}_0 \in [0,1]^3$ and a
neighborhood which includes ${\bf x}_0$, but is not necessarily centered at ${\bf x}_0$. A linear
model is to be fit over this region. The
slope in the $X_1$ direction $\beta_1$, is estimated using least squares, giving us
$\widehat \beta_1$. The $t$-statistic for $\widehat \beta_1$ follows easily.
One can now imagine repeating this for all possible neighborhoods, and finding that
neighborhood which results in the largest absolute $t$-statistic for $\widehat \beta_1$.
We then imagine doing this for all possible values of ${\bf x}_0$.
In practice, a reasonable subset of all possible neighborhoods will be chosen, and models
fit over these regions. Here, we lay down a grid of values for each of the covariates,
based on evenly spaced quantiles of the data. The number of grid points can be different
for each covariate; we choose to have a larger number for the covariate of interest
(we use 15) than for the other covariates (we use 5). Thus, we have an interior grid made
of up of $15\times 5 \times 5 = 375$ covariate vectors. A neighborhood is formed by choosing
two of these points and using them as the corners of the region. Thus, there are a total
of $375 \times 374/2 = 70125$ potential neighborhoods; some will be discarded due to sparsity
of data.

\begin{figure}[ht]
\begin{center}
\subfigure[Raw $t$-stat Plot]
{
\psfig{file=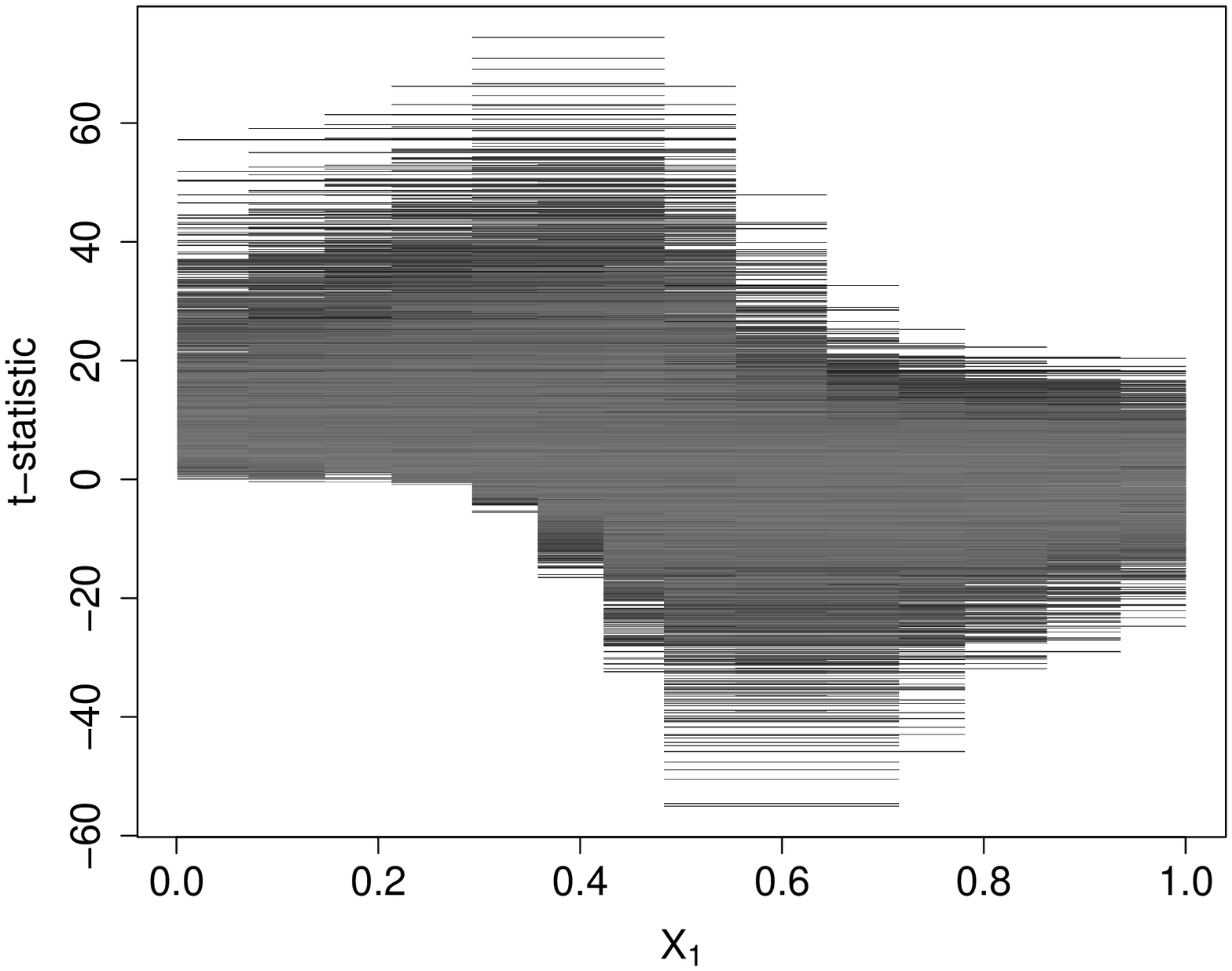,width=\twowidth}
\label{case1_rawplot}
}
\subfigure[$t$-stat Plot]
{
\psfig{file=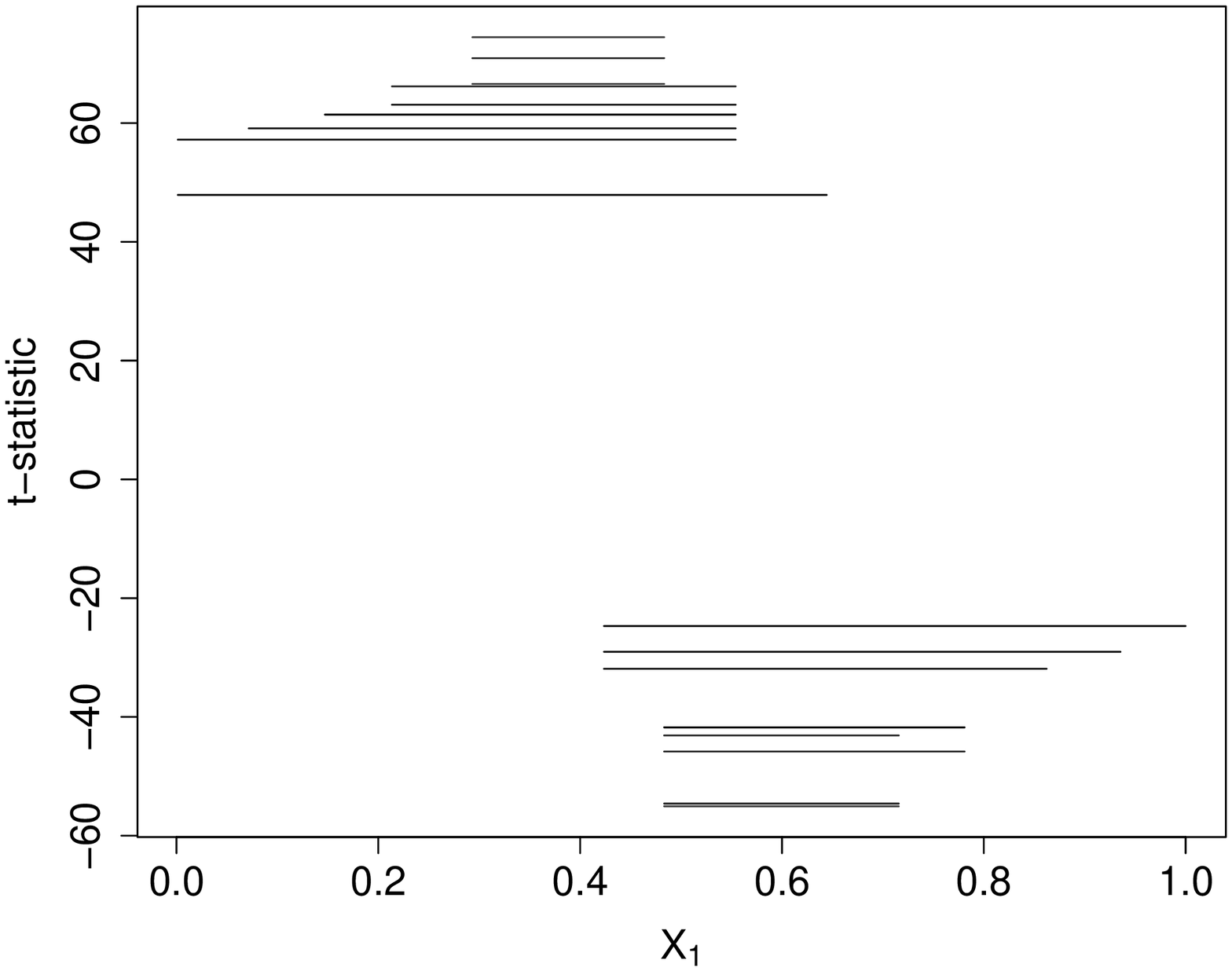,width=\twowidth}
\label{case1_tstatplot}
}
\caption{Figure \ref{case1_rawplot} the plot of $t$-statistics for each neighborhood fit for the simple case of Equation (\ref{simpmod}).
Figure \ref{case1_tstatplot} is the same, except it only shows 
those neighborhoods which maximize the absolute $t$-statistic for some covariate vector ${\bf x}_0$.}
\end{center}
\end{figure}

%We consider the basic example, Case 1, first.
We consider a simple case first. Set
\begin{equation}
   Y = f_1(X_1) + f_2(X_2) + f_3(X_3) + \epsilon
   \label{simpmod}
\end{equation}
where $\epsilon$ is normal with mean zero and variance $\sigma^2 = 0.02$. A sample of
size $n=1000$ is taken. The random variables $X_1, X_2, X_3$ are i.i.d. ${\cal U}(0,1)$.
Figure \ref{case1_rawplot} shows the raw results of the analysis. For each of the local linear
models, there is one line on the plot. The value on the vertical axis gives the $t$-statistic
for $\widehat \beta_3$. The range of values along the horizontal axis represents the range
of values of the variable of interest for that neighborhood. The shade of the line indicates the
proportion of the covariate space in the other two variables covered by that neighborhood. For
example, the darker lines indicate regions which, regardless of their width in the $X_1$ direction,
cover most of the $X_2$ and $X_3$ space.

%\begin{figure}[ht]
%\begin{center}
%\psfig{file=case1_rawplot.eps,width=\onecolwidth}
%\caption{The raw plot of $t$-statistics for each neighborhood fit for the simple case of Equation (\ref{simpmod}).}
%\label{case1_rawplot}
%\end{center}
%\end{figure}

%\begin{figure}[ht]
%\begin{center}
%\psfig{file=case1_tstatplot.eps,width=\onecolwidth}
%\caption{The $t$-statistic plot for the simple case, only showing those neighborhoods which maximize the absolute
%$t$-statistic for some covariate vector ${\bf x}_0$.}
%\label{case1_tstatplot}
%\end{center}
%\end{figure}

Figure \ref{case1_tstatplot} is the same as Figure \ref{case1_rawplot} 
except that the only regions plotted are
those that yield the maximum absolute $t$-statistic for some ${\bf x}_0$. 
A convenient computational
aspect of this approach is that wide ranges covariate values $x_0$ will share the same region: We
can describe fixing $x_0$ and finding the optimal region for that $x_0$, but in fact we only need
to fit all of these 70125 candidate models, calculate the $t$-statistic for each, and then discard
any fit which is not the maximum for some $x_0$. In this example there are 18 regions remaining
following this process.
The regions represented in this plot are the previously defined {\em features}, 
${\cal S}_1, {\cal S}_2, \ldots, {\cal S}_r$. We call this plot the ``$t$-statistic plot for
variable of interest $X_1$.''

\begin{figure}[ht]
\begin{center}
\psfig{file=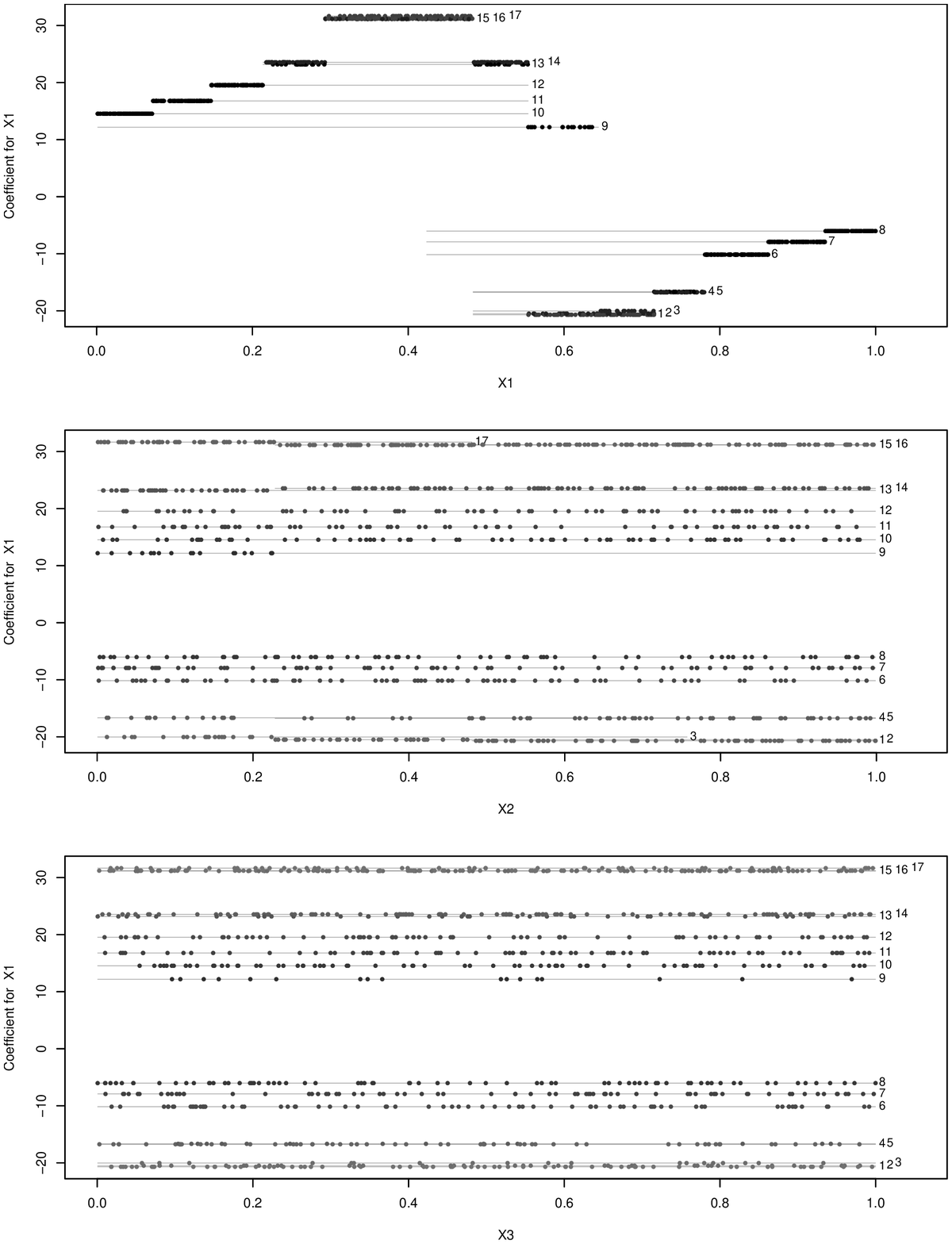,width=6.5in}
\caption{The feature plot for variable $X_1$, for the simple case.}
\label{case1_featureplot}
\end{center}
\end{figure}

A more useful way of plotting the remaining regions is shown in Figure \ref{case1_featureplot}. 
We will refer to this graph as the ``feature plot for variable of interest $X_1$.''
In each plot, there is one light, horizontal line for each feature; note that they
are labeled with numbers going from 1 to 18.  The vertical axis gives $\widehat \beta_1$ for
that local model. Each dot is an observed data point, and the shade of the dots on one line
(i.e., in one region) again represents the extent of that region in the other two covariates.

From this plot, one can pick out dependencies in the data. Consider the regions labeled ``1,''
``2,'' and ``3.'' For each of these, the range of values of $X_1$ in the region is approximately
0.5 to 0.7 (look at the third plot). 
This is capturing the steep downslope in $f_1(x)$ for $x$ in that range. But, the
dependence between $X_2$ and $Y$ (characterized here by the slope $\beta_2$) varies with $X_2$.
Ignoring this fact when modeling $Y$ as a function of $X_1$ would mask this downslope. This
is seen clearly when we look at these three regions in the second plot in Figure \ref{case1_featureplot}.
The three regions correspond to different ranges of values of $X_2$: Region ``1'' is approximately
0.5 to 0.75, where $f_2(x)$ is starting to level off, region ``2'' is extends up to 1.0, where $f_2(x)$
is almost flat, and region ``3'' goes from 0 up to 0.5, where $f_2(x)$ is the steepest. 

Since the
relationship between $X_3$ and $Y$ is linear, there is no such pattern to be found in the last of
the three plots. The light horizontal lines in this first plot show that our procedure chooses the
largest possible bandwidths in the direction $X_1$, that is, the response $Y$ is modeled to be linear
in this direction, as we know it should be.

\begin{figure}[ht]
\begin{center}
\subfigure[CV Plot]
{
\psfig{file=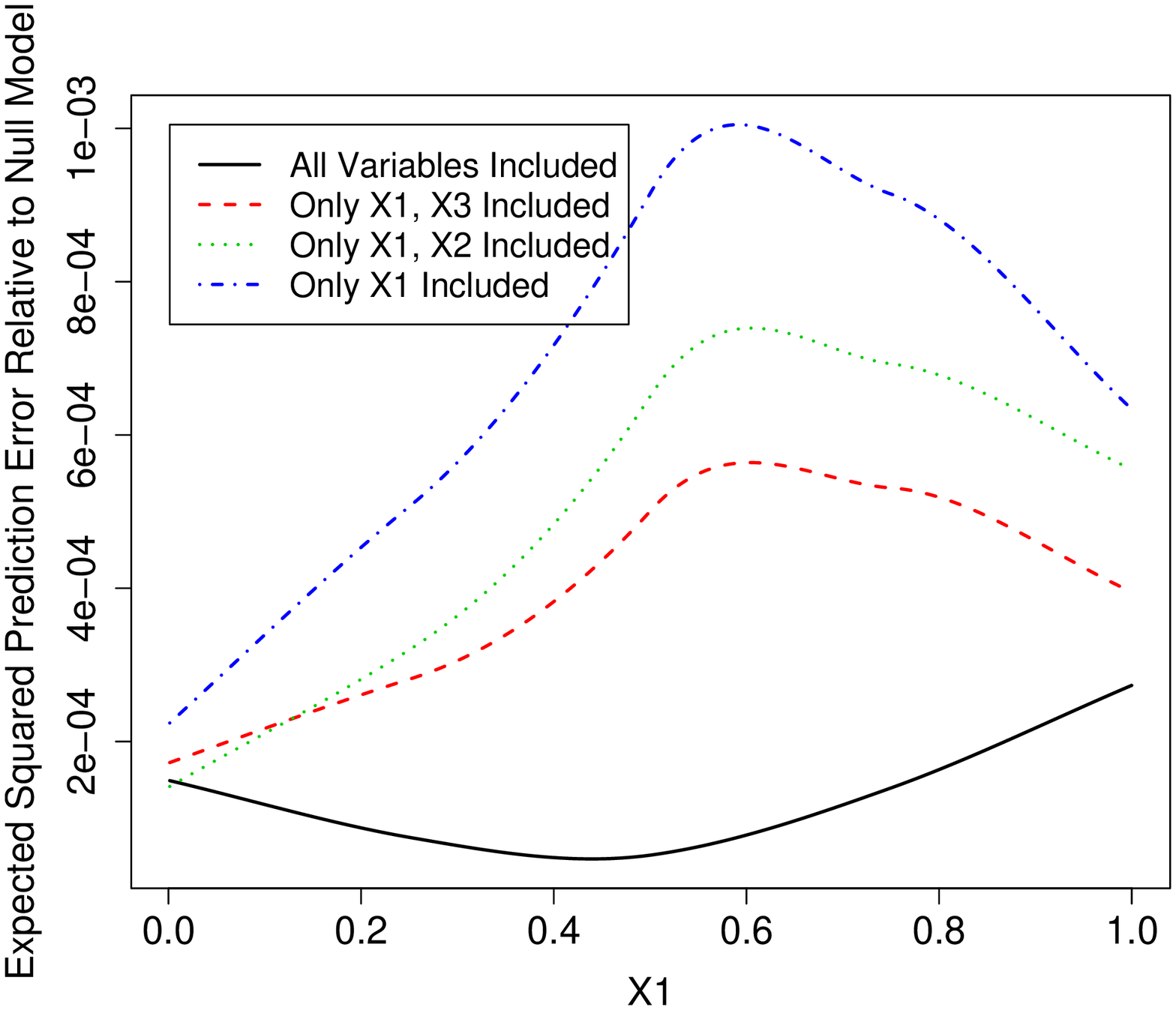,width=\twowidth}
\label{case1_cvplot}
}
\subfigure[Slope Plot]
{
\psfig{file=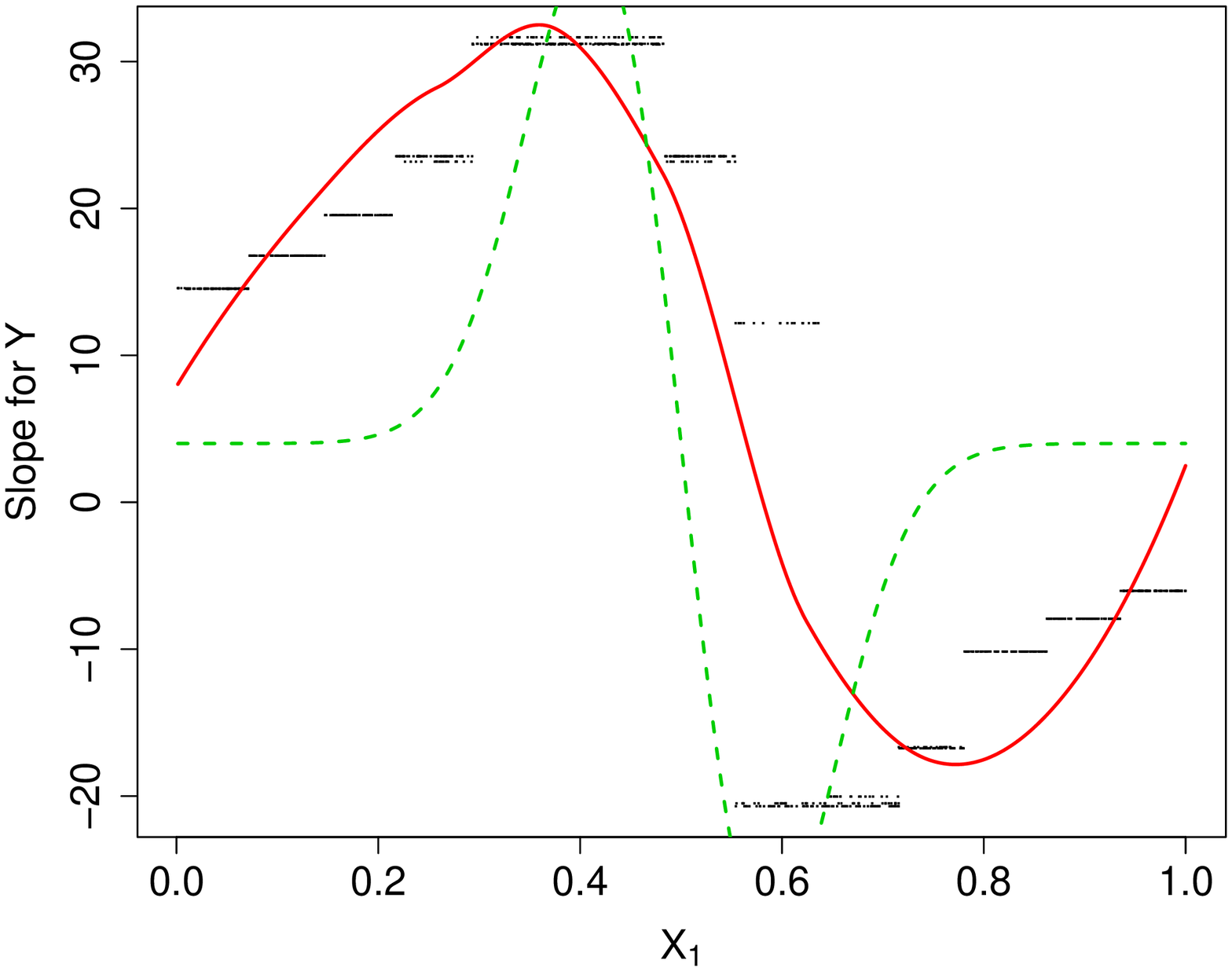,width=\twowidth}
\label{case1_slopeplot}
}
\caption{More results from the analysis of the simple model. Figure \ref{case1_cvplot} is 
the CV plot, comparing competing models for the purposes of variable selection.
Figure \ref{case1_slopeplot} is the slope plot for variable $X_1$.}
\end{center}
\end{figure}

Figure \ref{case1_cvplot} shows the estimates of the expected squared prediction error
$\gamma(\cdot,\widehat \mu_j)$ for each of four models, relative to this quantity for the ``null'' model,
the model which only uses the mean of $Y$ within that region to predict the response. 
We call this plot the ``CV plot.'' In this case, the best choice is to use all three predictors, as
evidenced by the solid black line being the lowest for all values of $X_1$.
Figure \ref{case1_slopeplot} (the ``slope plot'') 
shows, for each $X_1$ value in the data set, the estimated slope $\widehat \beta_1$
for the region ${\cal S}_k'$ in which that observation lies. In practice, if the variable selection procedure
chose a simpler model for a particular subset ${\cal S}_k$, that model would be used in the slope plot.

%\begin{figure}[ht]
%\begin{center}
%\psfig{file=case1_cvplot.eps,width=\onecolwidth}
%\caption{The CV plot for the simple case, comparing competing models for the purposes of variable selection.}
%\label{case1_cvplot}
%\end{center}
%\end{figure}

%\begin{figure}[ht]
%\begin{center}
%\psfig{file=case1_slopeplot.eps,width=\onecolwidth}
%\caption{The slope plot for variable $X_1$ for the simple case.}
%\label{case1_slopeplot}
%\end{center}
%\end{figure}

\subsection{Other Cases}

\begin{figure}
\begin{center}
\subfigure[CV Plot]
{
\psfig{file=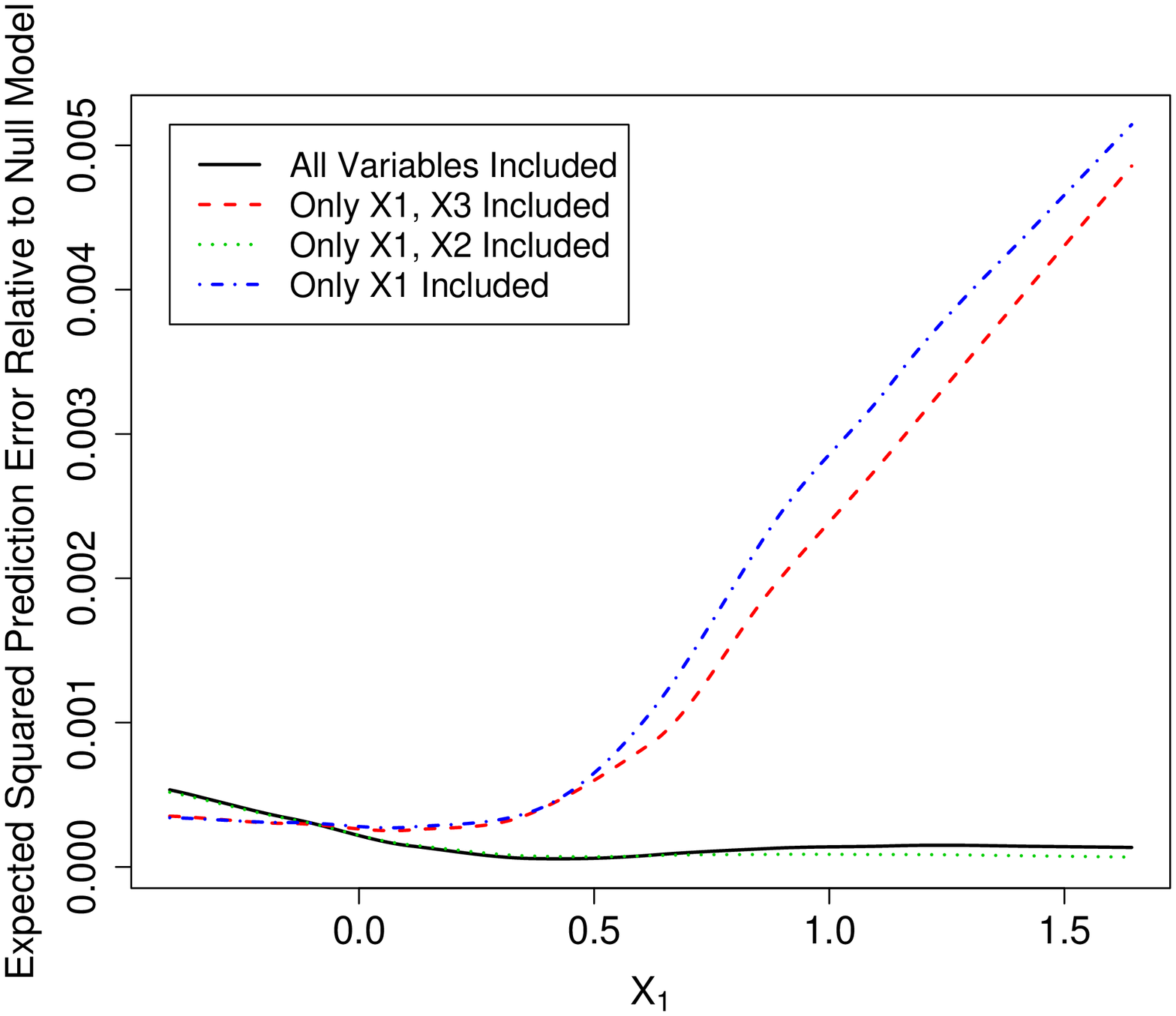,width=\twowidth}
\label{case6_cvplot}
}
\subfigure[CV Plot]
{
\psfig{file=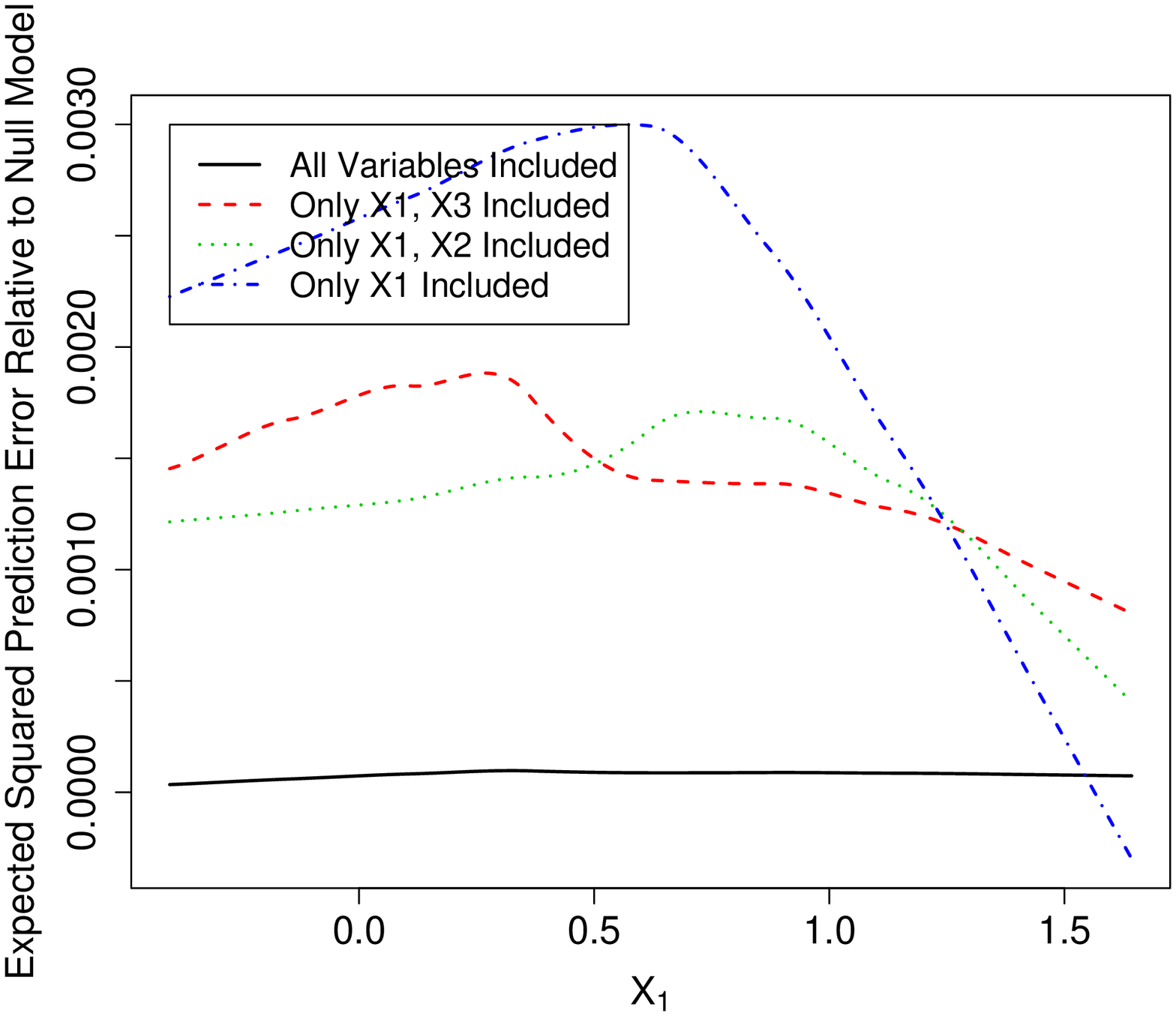,width=\twowidth}
\label{case7_cvplot}
}
\caption{Examples from analyses of extensions of the simple model. Figure \ref{case6_cvplot} is
the CV plot for the case where $X_1$ and $X_3$ are dependent, but $Y$ is not a function
of $X_3$, comparing competing models for the purposes of variable selection.
Figure \ref{case7_cvplot} is 
the CV plot for the case where $X_1$ and $X_3$ are dependent, but $Y$ is a function of $X_3$, not of $X_1$.}
\end{center}
\end{figure}

We considered alterations to the simple simulation model described
in the previous section. First, consider a case with $Y = f_1(X_1) + f_2(X_2) + \epsilon$,
so that $Y$ is not a function of $X_3$, but now take $(X_1,X_3)$ to be bivariate
normal, each with mean 0.5 and SD 1, and with correlation $\sqrt{0.5}$. $X_2$ is
still ${\cal U}(0,1)$, and independent of $X_1$ and $X_3$. Again, $n=1000$.
Figure \ref{case6_cvplot} shows the CV plot for this case. Note that since
conditional on $X_1$, $X_3$ and $Y$ are independent, the variable selection procedure
is indicating that $X_3$ could be excluded.
Contrast this with the second case where $(X_1, X_3)$ have the same bivariate
normal distribution, but now $Y=f_2(X_2) + f_3(X_3) + \epsilon$; see Figure \ref{case7_cvplot}
for the CV plot. Here, the best choice is to include all three variables since excluding
$X_3$ would lead to misleading conclusions regarding the strength of the relationship between
$X_1$ and $Y$.

%\begin{figure}[ht]
%\begin{center}
%\psfig{file=case6_cvplot.eps,width=\onecolwidth}
%\caption{The CV plot for the case where $X_1$ and $X_3$ are dependent, but $Y$ is not a function
%of $X_3$, comparing competing models for the purposes of variable selection.}
%\label{case6_cvplot}
%\end{center}
%\end{figure}

%\begin{figure}[ht]
%\begin{center}
%\psfig{file=case7_cvplot.eps,width=\onecolwidth}
%\caption{The CV plot for the case where $X_1$ and $X_3$ are dependent, but $Y$ is a function of $X_3$, not of $X_1$.}
%\label{case7_cvplot}
%\end{center}
%\end{figure}

\subsection{Analysis of Currency Exchange Data}

\begin{figure}
\begin{center}
\psfig{file=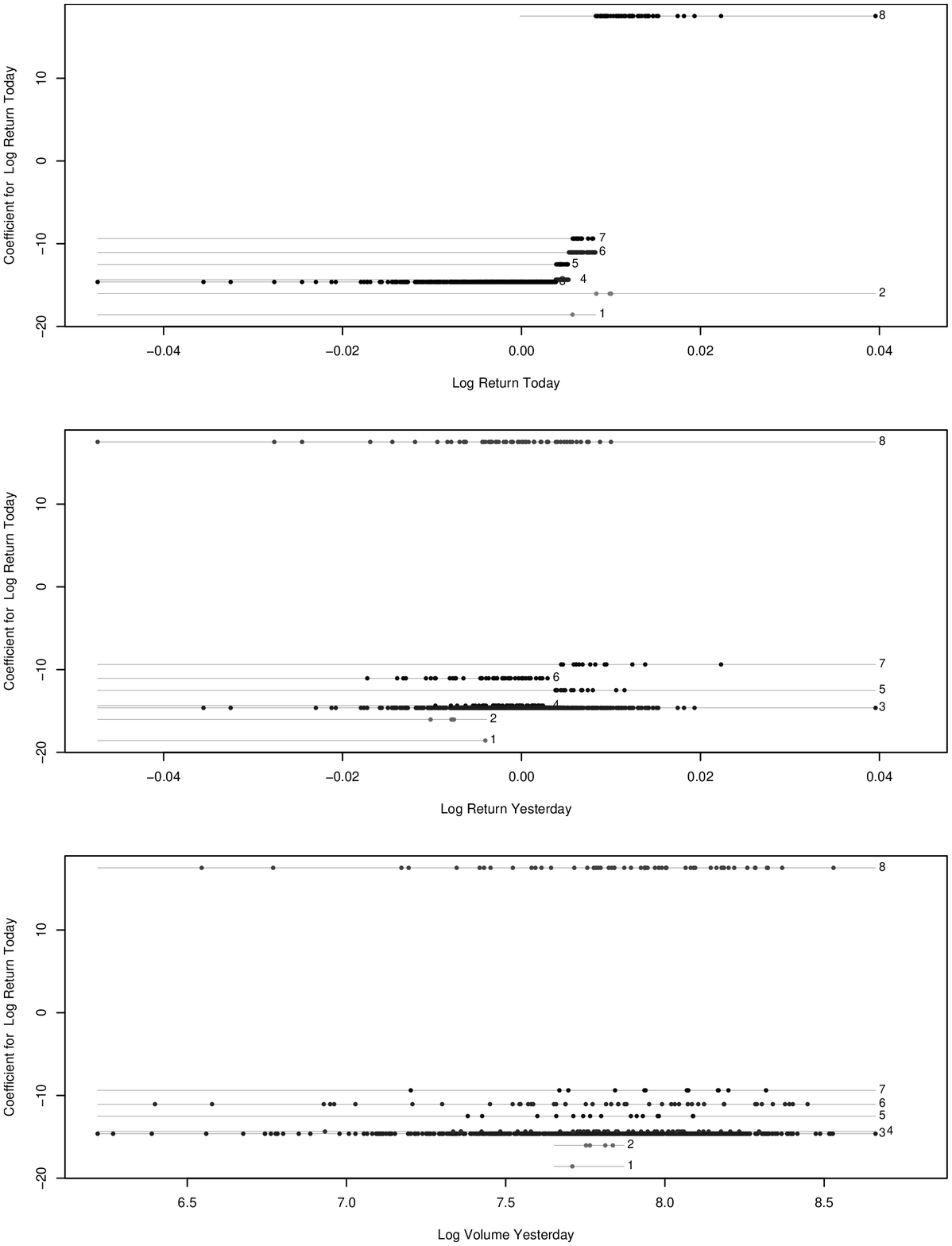,width=6.5in}
\caption{The feature plot for today's log volume for currency exchange data.}
\label{yenanal_featureplot}
\end{center}
\end{figure}

The original motivation for this study was to address a question in the
study of currency exchange rates regarding the nature of the relationship
between the volume of trading and the return.
Data were obtained on the Japanese Yen to U.S. dollar exchange rate for the
period of January 1, 1992 to April 14, 1995, a total of 1200 trading days.
The response variable used was today's log volume with three covariates:
1) today's log return, 2) yesterday's log volume, and 3) yesterday's log
return. The first of these, today's log return, is set as the covariate
of interest. Figure \ref{yenanal_featureplot} shows the feature plot
for this data set. One interesting result is that when today's log
return is positive, the coefficient for today's log return
is positive; see feature 8 at the top of the plot.
And when today's log return is negative, those coefficients are mostly negative;
see features 1,3,4, and so forth, on the left side of the plot.
This confirms a prediction of \citet*{Karpoff1987}.

\begin{figure}
\begin{center}
\subfigure[CV Plot]
{
\psfig{file=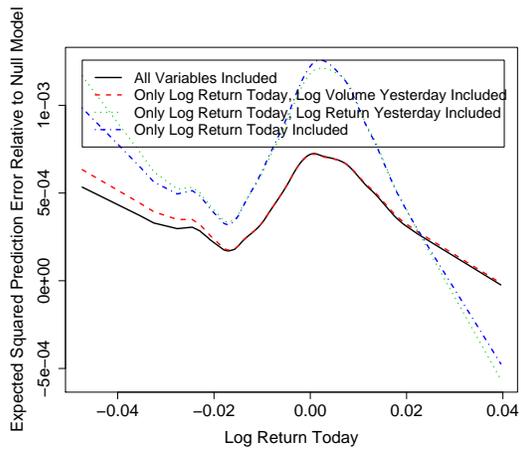,width=\twowidth}
\label{yenanal_cvplot}
}
\subfigure[Slope Plot]
{
\psfig{file=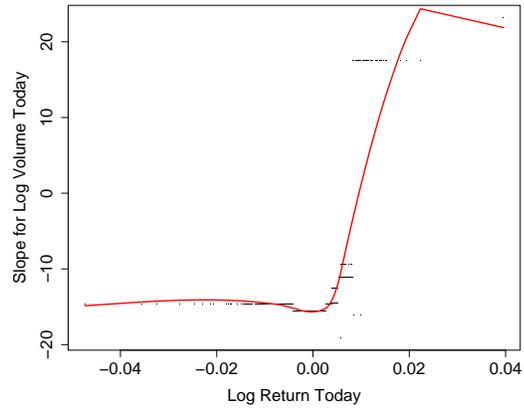,width=\twowidth}
\label{yenanal_slopeplot}
}
\subfigure[Level Plot]
{
\psfig{file=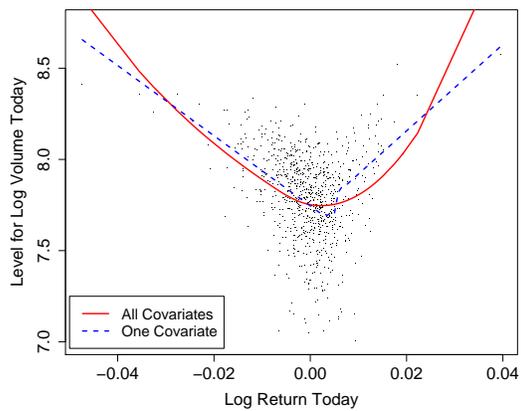,width=\twowidth}
\label{yenLVT_level}
}
\subfigure[Level Plot]
{
\psfig{file=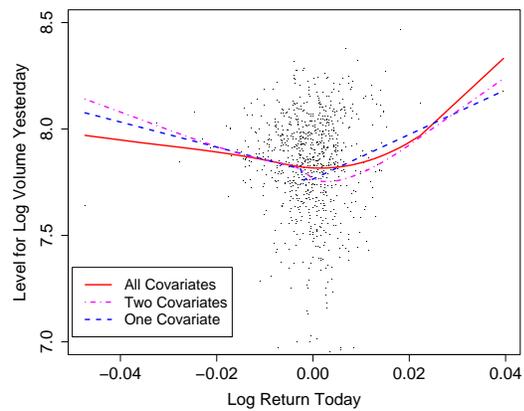,width=\twowidth}
\label{yenLVY_level}
}
\caption{Results from the analysis of the currency exchange data: 
The CV plot (Figure \ref{yenanal_cvplot}),
the slope plot (Figure \ref{yenanal_slopeplot}), and
the level plot (Figure \ref{yenLVT_level}) from the analysis
where the response is today's log volume.
Figure \ref{yenLVY_level} is the level plot for the 
using yesterday's log volume as the response.}
\end{center}
\end{figure}

When variable selection is applied, we see that there is evidence
that yesterday's log return is not important; see Figure \ref{yenanal_cvplot}.
Finally, the slope plot (Figure \ref{yenanal_slopeplot}) shows again how
the estimated slope $\widehat \beta_1$ for the coefficient for today's log return
abruptly switches from negative to positive once that variable becomes positive.
This is an important finding that may have been missed under a standard multiple
regression fit to this surface. 
This becomes clearer when inspecting the {\it level plot} shown in Figure \ref{yenLVT_level}. 
The plot shows one dot for each vector ${\bf x}$ in the data set, the horizontal 
axis gives today's log
return, and the vertical axis is the fitted value $\widehat \mu(x)$. The superimposed
solid curve is a smooth of the plot. 
For comparison, also plotted (as a dashed line) is the level curve 
for the $d=1$ case where the only covariate is today's log return, as in \citet*{Karpoff1987}.

Figure \ref{yenLVY_level} is the level plot from the same analysis,
except now making yesterday's log volume the response, and using today's log
volume as one of the explanatories. Once again, we see some evidence of the
``Karpoff effect.''
An interesting aspect of this plot is that when all three covariates are
included in the analysis, the slope on the left of the level curve
is significantly smaller than the
slope on the right. This seems to imply that there would be a way to
exploit yesterday's log volume in an effort to predict whether 
today's log return will be positive or negative.
But, this artifact is removed by excluding
today's log volume as one of the covariates:
See the level curve corresponding to ``Two Covariates.'' Using only
yesterday's information, it now does not seem possible to
utilize yesterday's log volume to predict if today's log return will be
positive or negative.

%\begin{figure}
%\begin{center}
%\psfig{file=franksanal_cvplot.eps,width=\onecolwidth}
%\caption{The CV plot for the analysis of the currency exchange data.}
%\label{franksanal_cvplot}
%\end{center}
%\end{figure}

%\begin{figure}
%\begin{center}
%\psfig{file=franksanal_slopeplot.eps,width=\onecolwidth}
%\caption{The slope plot for the analysis of the currency exchange data.}
%\label{franksanal_slopeplot}
%\end{center}
%\end{figure}

%\begin{figure}
%\begin{center}
%\psfig{file=franksanal_levelplot.eps,width=\onecolwidth}
%\caption{The level plot for the analysis of the currency exchange data.}
%\label{franksanal_levelplot}
%\end{center}
%\end{figure}

%\begin{figure}
%\begin{center}
%\psfig{file=franksanal_levelLVY.eps,width=\onecolwidth}
%\caption{The level plot for the analysis of the currency exchange data, now using yesterday's log
%volume as the response.}
%\label{franksanal_levelplotLVY}
%\end{center}
%\end{figure}

\section{Appendix}

\subsection{Proof of Lemma \ref{lemma1}}
The following implies Lemma \ref{lemma1}.

\begin{lem}\label{lemA1}
Let ${\rm RSS}^{\bf h} \equiv \sum^{\bf h} [Y_i - \widehat \mu_L({\bf X}_i)]^2$,
then
\begin{equation}
   E^{\bf h}\!\left({\rm RSS}^{\bf h}\right)
   = \left[n {\cal V}\!\left({\bf h}\right) - d\right]\sigma^2 + n {\cal V}\!\left({\bf h}\right)
   E^{\bf h}\!\left[\mu\!\left({\bf X}\right)-\mu_L\!\left({\bf X}\right)\right]^2.
\end{equation}
\end{lem}

\begin{proof}
If we condition on ${\mathbb X} \equiv ({\bf X}_1, \ldots, {\bf X}_n)$, then ${\bf X}_0 \in N_{\bf h}({\bf x}_0)$
no longer are random events. We can adapt \citet{Hastie1987}, equation (16) to find
\begin{equation}
   E^{\bf h}\!\left({\rm RSS}^{\bf h} | {\mathbb X}\right)
   = n_{\bf h} \sigma^2 + \Sigma^{\bf h}\left[ \mu\!\left({\bf X}_i\right) - E\!\left(\widehat \mu_L\!\left({\bf X}_i\right)\right)\right]^2
   - d \sigma^2.
\end{equation}
Here $E^{\bf h}(\widehat \mu_L({\bf X}_i)) = \mu_L({\bf X}_i)$ and Lemma \ref{lemA1} follows by the
iterated expectation theorem. Lemma \ref{lemma1} is a special case with $\mu$ and $\mu_L$ computed under
the null hypothesis.
\end{proof}

\subsection{Proof of Proposition \ref{prop5}}
\begin{proof}
We compute expected values using the joint density
\begin{equation}
   f^{\bf h}\!\left({\bf x}\right) = f\!\left({\bf x}\right)
   {\bf 1}\!\left[{\bf x} \in N_{\bf h}\!\left({\bf x}_0\right)\right]
   \bigg/ {\cal V}\!\left({\bf h}\right)
\end{equation}
of ${\bf X}$ given ${\bf X} \in N_{\bf h}({\bf x}_0)$. Thus, with $k=j$,
\begin{eqnarray}
   {\rm Cov}^{\bf h}\!\left(X_j, W_j\!\left(\frac{X_j - x_{0j}}{\theta_j}\right)\right)
   & \!\!=\!\! & \int \left[x_j - {\rm E}^{\bf h}\!\left(X_j\right)\right]
   W_j\!\left(\frac{x_j-x_{0j}}{\theta_j}\right) f^{\bf h}\!\left({\bf x}\right) d{\bf x} \nonumber \\
   & \!\!=\!\! & \int \left[x_j - {\rm E}^{\bf h}\!\left(X_j\right)\right]
   W_j\!\left(\frac{x_j-x_{0j}}{\theta_j}\right) f_j^{\bf h}\!\left(x_j\right) dx_j,
\end{eqnarray}
where
\begin{equation}
   f_j^{\bf h}\!\left(x_j\right) \equiv \int f^{\bf h}\!\left({\bf x}\right) d{\bf x}_{(-j)},
\end{equation}
with ${\bf x}_{(-j)} = (x_1,\ldots,x_{j-1},x_{j+1},\ldots,x_d)^T$, is the marginal
density of $X_j$ given ${\bf X} \in N_{\bf h}({\bf x}_0)$.
The change of variables $s_j = (x_j - x_{0j})/\theta_j$ gives
\begin{equation}
   {\rm Cov}^{\bf h}\!\left(X_j, W_j\!\left(\frac{X_j-x_{0j}}{\theta_j}\right)\right)
   = \theta_j \int_{-h_j}^{h_j} \left[{\bf x}_{0j}-E^{\bf h}\!\left(X_j\right)
   + s_j\right] W_j\!\left(s_j\right) f_j\!\left(x_{0j}+\theta_j s_j\right) ds_j.
\end{equation}
A similar argument for $k \neq j$ shows that $\gamma_k {\rm Cov}^{\bf h}(X_j, Y)=O(\max_k\{\gamma_k \theta_k\})$.
Now (a) follows because the terms in ${\rm EFF}_1$ other than $\gamma_k {\rm Cov}^{\bf h}(X_j,Y)$ are fixed as
$\gamma_k \theta_k \rightarrow 0$. Finally, (b) follows from Proposition \ref{consest}.
\end{proof}

\subsection{Proof of Theorem \ref{thm1}}
\begin{proof}
Because of the independence of the $X_1$ and $X_2, \ldots, X_d$,
%(a) follows from the proof of Theorem 2.1 in \citet{Doksum2006}
%for the $d=1$ case. Next note that
\begin{equation}
  n^{-1/2} {\rm EFF}_1\!\left({\bf h},{\bf x}_0\right) =
   \left[\sigma^2 + \sigma_{\bf L}^2\!\left(\hmone\right)\right]^{1/2}
  {\cal V}^{1/2}\!\left(h_1\right) {\cal V}^{1/2}\!\left(\hmone\right)
  \left[{\rm SD}^{h_1}\!\left(X_1\right)\right]^{-1}
  {\rm Cov}^{h_1}\!\left(X_1,Y\right)
\end{equation}
where ${\cal V}(h_1) \equiv P(x_{01}-h_1 \leq X_1 \leq x_{01} + h_1)$.
%Thus for $j \geq 2$, $h_j = \infty$ maximizes the efficacy, in which case
%\[
%   n^{-1/2} {\rm EFF}_1\!\left({\bf h},{\bf x}_0\right) \propto \sigma^{-1}
%   \gamma_1 \theta_1^2 h_1^{-3/2} m_1\!\left(W_1\right) f_1^{1/2}\!\left(x_{01}\right)
%\]
%by
The result now follows from the proof of Theorem 2.1 in \citet{Doksum2006}.
\end{proof}

\subsection{Proof of Theorem \ref{thm2}}
\begin{proof}
The proof can be constructed by using the fact that for small $|{\bf h}|$, $X_1,\ldots, X_d$ given
${\bf X} \in N_{\bf h}({\bf x}_0)$ are approximately independent with ${\cal U}(x_{0j}-h_j,\:x_{0j}+h_j)$,
$1 \leq j \leq d$, distributions. This can be seen by Taylor expanding $f({\bf x})$ around $f({\bf x}_0)$
and noting that
\begin{eqnarray}
   P\!\left({\bf X} \in A \:|\: {\bf X} \in N_{\bf h}\!\left({\bf x}_0\right)\right) & = &
   \int_{A \:\cap\: N_{\bf h}\!\left({\bf x}_0\right)}
   f\!\left({\bf x}\right) d{\bf x} \:\bigg/\! \int_{N_{\bf h}\!\left({\bf x}_0\right)}
   f\!\left({\bf x}\right) d{\bf x} 
%  + o\!\left(|{\bf h}|^2\right) 
   \nonumber \\
   & = &
   \int_{A \:\cap\: N_{\bf h}\!\left({\bf x}_0\right)}
   d{\bf x} \:\bigg/ \!\int_{N_{\bf h}\!\left({\bf x}_0\right)}
   d{\bf x} + o\!\left(|{\bf h}|^2\right).
\end{eqnarray}
A similar approximation applies to moments. Now use the proof of Theorem \ref{thm1} with appropriate small
error terms.
\end{proof}

\bibliography{kadcms}

\end{document}